\documentclass{article}

\usepackage{graphicx}
\usepackage{psfrag}
\usepackage{amsmath, amssymb, amsthm}
\usepackage{epstopdf}
\usepackage{enumerate}

\newtheorem{theorem}{Theorem}
\newtheorem{conjecture}[theorem]{Conjecture}
\newtheorem{proposition}[theorem]{Proposition}
\newtheorem{lemma}[theorem]{Lemma}
\newtheorem{corollary}[theorem]{Corollary}
\newtheorem{definition}{Definition}
\newtheorem{question}{Question}

\newcommand{\floor}[1]{\left\lfloor #1 \right\rfloor}

\newcommand{\reals}{\mathbb{R}}
\renewcommand{\L}{\mathcal{L}}

\title{Sampling Lissajous and Fourier Knots}
\author{Adam Boocher\footnote{Supported by NSF REU grant DMS-0453284.}\\
University of Notre Dame\\
\\
 Jay Daigle$^*$\\
 Pomona College\\
 \\
  Jim Hoste$^*$\\
  Pitzer College\\
  \\
   Wenjing Zheng$^*$\\
   University of California,  Berkeley}
\date{\today}
\begin{document}
\maketitle
\begin{abstract}
A {\it Lissajous knot} is one that can be parameterized as $$K(t)=\left (\cos (n_x t+\phi_x), \cos (n_y t+\phi_y), \cos (n_z t+\phi_z)\right )$$ where the {\it frequencies} $n_x, n_y$, and $n_z$ are relatively prime integers and the {\it phase shifts} $\phi_x, \phi_y$ and $\phi_z$ are real numbers. Lissajous knots are highly symmetric, and for this reason, not all knots are Lissajous. We prove several theorems which allow us to place bounds on the number of Lissajous knot types with given frequencies and  to efficiently sample all possible Lissajous knots with a given set of frequencies. In particular, we systematically tabulate all Lissajous knots with  small frequencies and as a result substantially enlarge the tables of known Lissajous knots.

A {\it Fourier-$(i, j, k)$ knot} is similar to a Lissajous knot except that the $x, y$ and $z$ coordinates are now each described by a sum of $i, j$ and $k$ cosine functions respectively. According to Lamm, every  knot is a Fourier-$(1,1,k)$ knot for some $k$. By randomly searching the set of Fourier-$(1,1,2)$ knots we find that all 2-bridge knots up to 14 crossings are either Lissajous or Fourier-$(1,1,2)$ knots. We show that all twist knots are Fourier-$(1,1,2)$ knots and give evidence suggesting that all torus knots are Fourier-$(1,1,2)$ knots.

As a result of our computer search, several knots with relatively small crossing numbers are identified as potential counterexamples to interesting conjectures.
\end{abstract}
\pagebreak
\section{Introduction}
A {\it Lissajous} knot $K$ in $\reals^3$ is a knot that has a
parameterization $K(t)=(x(t), y(t), z(t))$ given by
$$
     \begin{array}{rcl}
     x(t) & = & \cos(n_xt + \phi_x) \\
     y(t) & = & \cos(n_yt + \phi_y) \\
     z(t) & = & \cos(n_zt + \phi_z) \\
     \end{array}
$$
where $0 \leq t \leq 2\pi$, $n_x, n_y$, and $n_z$ are integers, and
$\phi_x, \phi_y, \phi_z \in \reals$. \label{lissajous definition}

Lissajous knots were first studied in \cite{BHJ1994} where some of their elementary properties were established. Most notably, Lissajous knots enjoy a high degree of symmetry. In particular, if the three {\it frequencies} $n_x, n_y$ and $n_z$ (which must be pairwise relatively prime---see \cite{BHJ1994}) are all odd, then the knot is strongly plus amphicheiral. If one of the frequencies is even, then the knot is 2-periodic, with the additional property that it links its axis of rotation once. These symmetry properties imply (strictly) weaker properties such as the fact that the Alexander polynomial of a Lissajous knot must be a square mod 2, which in turn implies that its Arf invariant must be zero. See  \cite{BHJ1994}, \cite{HZ2006} and \cite{Lamm1996} for details. Thus for example, the trefoil and figure eight knots are not Lissajous since their  Arf invariants are one. In fact, ``most'' knots are not Lissajous.

To date it is unknown if every knot which is strongly plus amphicheiral or 2-periodic (and links its axis of rotation once) is Lissajous. Several  knots with relatively few crossings exist which meet these symmetry requirements and yet are still unknown to be Lissajous or not. For example, according to \cite{HTW:1998} there are only three prime knots with 12 or less crossings which are strongly plus amphicheiral: 10a103 ($10_{99})$, 10a121 ($10_{123}$), and 12a427.   Here we have given knot names  in both the  Dowker-Thistlethwaite ordering of the Hoste-Thistlethwaite-Weeks table  \cite{HTW:1998} and, in parenthesis,  the Rolfsen~\cite{Rolfsen} ordering (for knots with 10 or less crossings). Symmetries of the knots in the Hoste-Thistlethwaite-Weeks table were computed using {\it SnapPea} as described in \cite{HTW:1998}. Of these three knots, only 10a103 ($10_{99}$) was previously reported as  Lissajous. (See \cite{Lamm} and \cite{Lamm1996}.) However we find 12a427 to be Lissajous.  (See Section~\ref{sampling} of this paper.) This leaves open the case of 10a121.  As a further example, there are exactly four 8-crossings knots which are 2-bridge, 2-periodic, and link their axis of rotation once. Despite our extensive searching (see Section~\ref{sampling}) only one of these knots turned up as Lissajous (and it had already been reported as such in \cite{Lamm}). Whether the other three are Lissajous remains unknown.

Lissajous knots are a subset of the more general class of {\it Fourier} knots. A Fourier-$(i,j,k)$ knot is one that can be parameterized as 
$$
\begin{array}{rcl}
     x(t) & = & A_{x,1}\cos(n_{x,1}t+\phi_{x,1})+...+A_{x,i}\cos(n_{x,i}t+\phi_{x,i}) \\
     y(t) & = & A_{y,1}\cos(n_{y,1}t+\phi_{y,1})+...+A_{y,j}\cos(n_{y,j}t+\phi_{y,j}) \\
     z(t) & = & A_{z,1}\cos(n_{z,1}t+\phi_{z,1})+...+A_{z,k}\cos(n_{z,k}t+\phi_{z,k}). \\
\end{array}
$$
Because any function can be closely approximated by a sum of cosines, every knot is a Fourier knot for some $(i,j,k)$. But a remarkable theorem of Lamm \cite{Lamm} states that in fact every knot is a Fourier-$(1,1,k)$ knot for some $k$. While $k$ cannot equal one for all knots (these are the Lissajous knots, and not all knots are Lissajous) could $k$ possibly be less than some universal bound $M$ for all knots? This seems unlikely, with the more reasonable outcome being that $k$ depends on the specific knot $K$. Yet no one has found a knot for which $k$ must be bigger than two!

If $K$ is a Fourier-$(1,1,k)$ knot  then its bridge number is less than or equal to the minimum of $n_x$ and $n_y$. (The bridge number of a knot $K$ can be defined as the smallest number of extrema on $K$ with respect to a given direction in $\mathbb R^3$, taken over all representations of $K$ and with respect to all directions. See \cite{BZ:2003} or \cite{Rolfsen} for more details.) Moreover, Lamm's proof is constructive and explicitly shows that if $K$ has bridge number $b$, then $K$ is a Fourier-$(1,1,k)$ knot  for some $k$ and with $n_x=b$. This raises several interesting questions. For any knot $K$, when expressed as a Fourier-$(1,1,k)$ knot, can   the minimum values of $n_x$ and $k$ be simultaneously realized? In particular, can a knot which is Lissajous and with bridge index $b$ be realized as a Lissajous knot with $n_x=b$? 

Let $\L(n_x,n_y,n_z)$ be the set of all Lissajous knots  with
frequencies $n_x, n_y, n_z$. (Throughout this paper we consider a knot and its mirror image to be equivalent.) One of the main goals of this paper is to investigate the set $\L(n_x,n_y,n_z)$. By a simple change of variables, $t \to t+c$, we may  alter the phase shifts. Therefore we will assume that $\phi_x=0$ in all that follows. This leaves the pair of parameters $(\phi_y, \phi_z)$ which vary within the {\it phase torus} $[0,2 \pi]\times [0,2 \pi]$. In Section~\ref{lissajous phase torus} we examine the phase torus and identify a finite number of regions in which the phase shifts must lie, with each region corresponding to a single knot type. We further show that a periodic pattern of knot types are produced as one traverses the phase torus. This allows us to prove
\begin{theorem} \label{first bound} Let $|\L(n_x, n_y, n_z)|$ be the number of distinct Lissajous knots with frequencies $(n_x,n_y,n_z)$. Then
$$|\L(n_x, n_y, n_z)|\le 2n_xn_y.$$
If furthermore $n_x=2$, then 
$$|\L(2, n_y, n_z)|\le 2n_y+1.$$
\end{theorem}

 There is also a periodicity that exists across frequencies and in Section~\ref{lissajous phase torus} we also prove
\begin{theorem} \label{periodicity in n_z}
$\L(n_x,n_y,n_z)\subseteq \L(n_x,n_y,n_z+2n_xn_y)$, with equality if
$n_z\ge 2 n_xn_y-n_y$.
\end{theorem}

Our analysis of the phase torus, together with these theorems 
 allow us to efficiently sample (with the aid of a computer) all possible Lissajous knots having two of the three frequencies bounded. Even with relatively small frequencies, the three natural projections of a Lissajous knot into the three coordinate planes can have a large number of crossings. (The projection into the $xy$-plane has $2n_xn_y-n_x-n_y$ crossings.) With frequencies of 10 or more, diagrams with hundreds of crossings result and many, if not most, knot invariants are computationally out of reach. Thus it becomes extremely difficult to compare different Lissajous knots with large frequencies, or to try to locate them in existing knot tables. However, if one frequency is two, the knot is 2-bridge and even with hundreds of crossings it is relatively simple to compute the identifying fraction $p/q$ by which 2-bridge knots are classified.

In Section~\ref{lissajous phase torus} we recall basic facts about Lissajous knots  and prove several theorems, including the two already mentioned, that will allow us to efficiently sample all Lissajous knots with two given frequencies. In Section~\ref{2-bridge knots} we then recall some basic facts about 2-bridge knots. Using these results we then report in Section~\ref{sampling} on our computer experiments. Theorems similar to those given in Section~\ref{lissajous phase torus} but for Fourier-$(1,1,k)$ knots would necessarily be much more complicated and we only begin the analysis of the phase torus for Fourier-$(1,1,2)$ knots in Section~\ref{fourier phase torus}.  Without the analogous results, we have not been able to rigorously sample Fourier knots. Instead, we have proceeded by two methods, either random sampling, or a sampling based on first forming a bitmap image of the phase torus and its singular curves.
However, even without exhaustive sampling,  our data show that all 2-bridge knots up to 14 crossings are Fourier-$(1,1,k)$ knots with $k\le 2$ and with $n_x=2$.

This research was carried out at the Claremont College's REU program in the summer of 2006. The authors thank the National Science Foundation and the Claremont Colleges for their generous support.

\section{The Phase Torus---Lissajous Knots}
\label{lissajous phase torus}

Suppose $K(t)$ is a Lissajous knot and consider its diagram in the $xy$-plane.
 Each crossing in this diagram corresponds to a double point in the $xy$-projection
given by a pair of times $(t_1, t_2)$, where $x(t_1) = x(t_2)$ and
$y(t_1) = y(t_2)$. The following lemma is given in   \cite{HZ2006}.
\begin{lemma}
Let $K(t)$ be a Lissajous knot. There are two types of time pairs
$(t_1, t_2)$ that give double points in the $xy$-projection:
$$ \begin{array}{rcl}
    \mbox{Type I}: & (t_1, t_2) = & \left( (-\frac{k}{n_x} + \frac{j}{n_y}) \pi - \frac{\phi_y}{n_y}, (\frac{k}{n_x} + \frac{j}{n_y}) \pi - \frac{\phi_y}{n_y} \right)
    \\
    &&\\
    & 1 \leq k \leq n_x - 1, & 1+\lfloor \frac{n_y}{n_x} k + \frac{\phi_y}{\pi}\rfloor \leq j \leq \lfloor 2 n_y - \frac{n_y}{n_x} k + \frac{\phi_y}{\pi} \rfloor \\
    & & \\

    \mbox{Type II}: & (t_1, t_2) = & \left( (-\frac{k}{n_y} + \frac{j}{n_x}) \pi - \frac{\phi_x}{n_x}, (\frac{k}{n_y} + \frac{j}{n_x}) \pi - \frac{\phi_x}{n_x}\right)\\
&&\\
    & 1 \leq k \leq n_y - 1, & 1+\lfloor \frac{n_x}{n_y} k + \frac{\phi_x}{\pi} \rfloor \leq j \leq \lfloor 2 n_x - \frac{n_x}{n_y} k + \frac{\phi_x}{\pi} \rfloor \\
\end{array}
$$
There are $n_x n_y - n_y$ double points of Type I, and $n_x n_y -
n_x$ double points of Type II.
\end{lemma}

\begin{figure}[t]
    \begin{center}
    \leavevmode
    \scalebox{.50}{\includegraphics{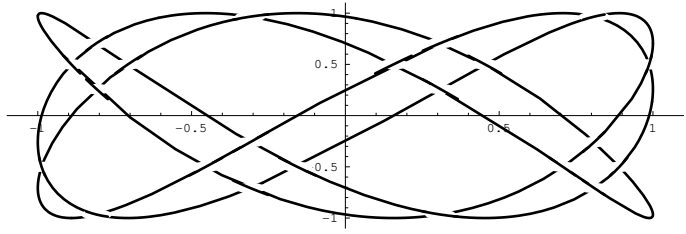}}
    \end{center}
\caption{A Lissajous knot with frequencies $(3,5,7)$ and corresponding phase shifts $(0,\pi/4, \pi/12)$. The Type I crossings appear in  two rows with five crossings each and the Type II crossings appear in four columns with three crossings each.}
\label{357knot}
\end{figure}

Figure~\ref{357knot} shows a Lissajous knot with frequencies $(3,5,7)$ and corresponding phase shifts $(0,\pi/4, \pi/12)$. Since all the frequencies are odd, this knot is symmetric through the origin. It is not hard to show that in general, the Type I  crossings line up in sets of size $n_y$ on $n_x-1$ horizontal lines, while the Type II crossings line up in sets of size $n_x$ on $n_y-1$ vertical lines. If $n_x=2$ there is a single row of Type I crossings, all of which lie on the $x$-axis and $n_y-1$ columns of Type II crossings with each column consisting of two crossings.

Not all phase shift pairs will generate curves that are knots. Assuming 
$\phi_x=0$, the knot $K(t)$ will intersect itself, and thus fail to be a knot, exactly
when the phase shifts satisfy
\begin{eqnarray}
\phi_z&=&\frac{n_z}{n_y}\phi_y + l\frac{\pi}{n_y}\mbox{, or }\label{Type I singular lines}\\
\phi_z&=&l \frac{\pi}{n_x}\mbox{, or }\label{Type II singular lines}\\
\phi_y&=&l\frac{\pi}{n_x}\label{vertical singular lines}
\end{eqnarray}
for some integer $l$.  Crossings of Type I become singular precisely when
Equation~\ref{Type I singular lines} holds; crossings of Type II
when Equaton~\ref{Type II singular lines} holds.  When Equation~\ref{vertical singular lines} holds, the entire $xy$-projection degenerates to an arc.  While this alone does not imply the knot has  points of self-intersection, this is indeed the case. See \cite{JP1998,BHJ1994,HZ2006} for more details. In Proposition~\ref{crossing diagonal}, we specifically identify which crossings become singular as the phase shifts move across these lines.

The slanted, horizontal and vertical lines given in Equations~\ref{Type I singular lines}--\ref{vertical  singular lines} obviously divide the phase torus into regions with each region defining one knot type. Thus there are only a finite number of knots types possible for a given set of frequencies. There is, however, a great deal of repetition in knot types as one traverses the phase torus due to the periodicity of the cosine function. The following theorem describes a nice choice of  ``fundamental domain'' on the phase torus to which we may restrict our attention.
\begin{theorem}
\label{fundamental domain}
Any knot in $\mathcal{L}(n_x,n_y,n_z)$ can be represented with $\phi_x=0$ and using some phase shift pair $(\phi_y, \phi_z)$ in $[0,\frac{\pi}{n_x}]\times[0,\pi]$.
\end{theorem}
\begin{proof}Define an equivalence relation $\sim$ on the phase torus for $\mathcal{L}(n_x,n_y,n_z)$ by $(\phi_y,\phi_z)\sim (\phi_y^\prime,\phi_z^\prime)$ if the Lissajous knot with phase shifts $(0,\phi_y,\phi_z)$ is the same as  the knot with phase shifts $(0,\phi_y^\prime,\phi_z^\prime)$, or its mirror image. Clearly
\begin{equation}
\label{pi shift}
(\phi_y,\phi_z)\sim (\phi_y \pm \pi,\phi_z)\sim (\phi_y,\phi_z \pm \pi).
\end{equation}
If $K\in \L(n_x,n_y,n_z)$, a change of variable $t\rightarrow t+\frac{\pi}{n_x}$ shows that $K$
is also parameterized as {$$
     \begin{array}{rcl}
     x & = & -\cos(n_xt) \\
     y & = & \cos(n_yt+\phi_y + \frac{n_y\pi}{n_x}) \\
     z & = & \cos(n_zt+\phi_z + \frac{n_z\pi}{n_x}). \\
     \end{array}
     $$}Therefore we also have 
\begin{equation}
\label{other shift}
(\phi_y,\phi_z)\sim (\phi_y +\frac{n_y\pi}{n_x}, \phi_z + \frac{n_z\pi}{n_x}).
\end{equation} 
Now since $n_x$
and $n_y$ are relatively prime, there are integers $k$ and $l$ with
${0\leq\phi_y+\frac{kn_y\pi}{n_x}-l\pi<\frac{\pi}{n_x}}$.  Repeatedly using (\ref{pi shift}) and (\ref{other shift}) we obtain 
$$(\phi_y,\phi_z)\sim(\phi_y+\frac{kn_y\pi}{n_x}-l\pi,\phi_z+\frac{kn_z\pi}{n_x})$$
The first coordinate is already in $[0,\frac{\pi}{n_x}]$; we can
shift the second coordinate by multiples of $\pi$ until it is in
$[0,\pi]$.  Thus an arbitrary point
$(\phi_y,\phi_z)$ is equivalent to some point in
$[0,\frac{\pi}{n_x}]\times[0,\pi]$, as desired.
\end{proof}

\begin{figure}[b!]

\psfrag{b}{\large $9/2$}
\psfrag{c}{\large $9/2$}
\psfrag{e}{\large $7/2$}
\psfrag{i}{\large $7/2$}
\psfrag{j}{\Large $\pi$}
\psfrag{k}{\Large $\frac{\pi}{2}$}
\psfrag{l}{\Large $0$}
\psfrag{m}{\Large $\frac{\pi}{2}$}
\psfrag{n}{\Large $\phi_z$}
\psfrag{a}{\Large {$\phi_y$}}
    \begin{center}
    \leavevmode
    \scalebox{.58}{\includegraphics{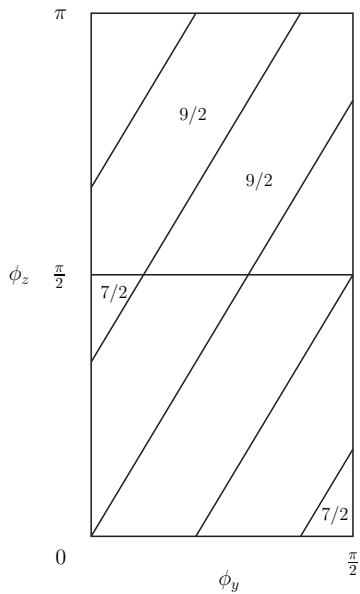}}
    \end{center}
\caption{The fundamental domain of the phase torus for $\L(2,3,5)$. Each region defines a single 2-bridge  knot which is identified by its classifying fraction $p/q$. Unlabeled regions define unknots.}
\label{235phasetorus}
\end{figure}

Figure~\ref{235phasetorus} shows the fundamental domain on the phase torus for $\L(2,3,5)$. The singular lines divide the domain into regions with each region determining a single knot type. Since these knots are all 2-bridge, we identify each with its classifying fraction $p/q$.

\pagebreak

Our next result specifically describes what happens as we cross a singular line of the type given in Equation~\ref{Type I singular lines} or \ref{Type II singular lines}.

\begin{proposition}
Let $K$ and $K'$ be two Lissajous knots with frequencies $(n_x, n_y, n_z)$ and phase shifts $(\phi_y,\phi_z)$
and $(\phi_y',\phi_z')$ respectively.
\begin{enumerate}
\item Suppose   $(\phi_y,\phi_z)$
and $(\phi_y',\phi_z')$ lie in two adjacent 
regions of the phase torus separated by a diagonal line $L$ given by
$\phi_z=\frac{n_z}{n_y}\phi_y+l\frac{\pi}{n_y}$.  Then $K$ and $K'$ differ by changing all  Type~I crossings with parameters $(k,j)$ such that with $jn_z+l\equiv0 \mod n_y$. The number of such crossings is $n_x-1$.
\item Suppose   $(\phi_y,\phi_z)$
and $(\phi_y',\phi_z')$ lie in two adjacent 
regions of the phase torus separated by a horizontal line $L$ given by $\phi_z=\l\frac{\pi}{n_x}$.  Then $K$ and $K'$ differ by changing all  Type~II crossings with parameters $(k,j)$ such that $jn_z+l\equiv0 \mod n_x$. The number of such crossings is $n_y-1$.
\end{enumerate}
\label{crossing diagonal}
\end{proposition}
\begin{proof}
 If $(t_1, t_2)$ is a Type I crossing with parameters $(k,j)$, then   $z(t_1)=z(t_2)$ if and only if 
$$\cos(n_zt_1+\phi_z)=\cos(n_zt_2+\phi_z)$$
which will occur if and only if
\begin{equation}n_z(t_1-t_2)=2m\pi~~~~\mbox{or}~~~~n_z(t_1+t_2)+2\phi_z=2m'\pi\label{singular or}\end{equation}
for some integers $m,m'$.  For Type I crossings,
$$t_1-t_2=-\frac{2k}{n_x}\pi~~~~\mbox{and}~~~~
t_1+t_2=\frac{2j}{n_y}\pi-\frac{2\phi_y}{n_y}.$$
If (\ref{singular
or}) is to hold, then in the first case, we have
$$-\frac{2kn_z}{n_x}\pi=2m\pi$$
which is equivalent to $-kn_z=mn_x$. This is impossible since $n_x$ and $n_z$ are relatively prime and $1\le k \le n_x-1$.

In the second case, we have
$$n_z\left(2\frac{j}{n_y}\pi-2\frac{\phi_y}{n_y}\right)+2\phi_z=2m'\pi$$
which is equivalent to
$$\phi_z=\frac{n_z}{n_y}\phi_y+(m'n_y-jn_z)\frac{\pi}{n_y}.$$
Thus Type I crossings only become singular on lines of the form given in Equation~\ref{Type I singular lines} with $l=mn_y-jn_z$.

If $\phi_z=\frac{n_z}{n_x}\phi_y+l \frac{\pi}{n_y}+\varepsilon$ and $jn_z+l=mn_y$ for some integer $m$ then it is straightforward to check that
$$z(t_1)-z(t_2)=(-1)^{m}2 \sin \varepsilon \sin\frac{kn_z \pi}{n_x}.$$
Hence, as we move across the line $L$ by letting $\varepsilon$ go from a small positive value to a small negative value,  the difference $z(t_1)-z(t_2)$ changes sign. Thus the Type I crossings with parameters $(k,j)$ actually change from over to under or vice versa, rather than simply becoming singular and then ``rebounding'' to their original positions.

Finally, note that once $l$ is fixed this does not necessarily uniquely determine $j$ and thus the corresponding Type I crossing. If both $j n_z+l\equiv 0$ mod $n_y$ and $j'n_z+l\equiv 0$ mod $n_y$ then $j\equiv j'$ mod $n_y$ since $n_y$ and $n_z$ are relatively prime. If $n_x=2$, then $k=1$ and $j$ lies in an interval of length $n_y$.  Thus with $n_x=2$ we have that $j$ is uniquely determined by $l$ and a single crossing changes as we move across $L$. But if $n_x>2$  and $k=1$ then $j$ lies in an interval of length greater than $n_y$. Hence two admissible values, $j$ and $j+n_y$, are possible. Using $j$, we must have $1\le k\le \floor{\frac{n_x}{n_y}(j-\frac{\phi_y}{\pi})}$ while for $j+n_y$ we must have $1\le k\le  \floor{\frac{n_x}{n_y}(-j+n_y+\frac{\phi_y}{\pi})}$. Thus the total number of possible points $(k,j)$ is $ \floor{\frac{n_x}{n_y}(j-\frac{\phi_y}{\pi})}+\floor{\frac{n_x}{n_y}(-j+n_y+\frac{\phi_y}{\pi})}=n_x-1$.

A similar discussion handles the Type II crossings.
\end{proof}
\begin{corollary} Suppose $K$ and $K'$ are Lissajous knots with frequencies $(n_x, n_y, n_z)$ and phase shifts which belong to regions separated by $2 n_y$ singular lines of the type given in Equation~\ref{Type I singular lines}. Then all Type I crossings are the same for both knots.
\label{periodicity in box}
\end{corollary}
\begin{proof}
From Proposition~\ref{crossing diagonal} we know that crossing the line $\phi_z=\frac{n_z}{n_y}\phi_y+l\frac{\pi}{n_y}$ changes exactly those Type I crossings with parameters $(k,j)$ for which $j n_z+l\equiv0$ mod $n_y$. Thus if we cross the singular line corresponding to $l$ and then later cross the line corresponding to $l+n_y$ the same set of Type I crossings will first be changed and then changed back again. Hence, after crossing over $2n_y$ such lines all Type I crossings will be restored to their original position.
\end{proof}
If $n_x=2$ there is even more repetition due to additional symmetry as is shown in the following result.
\begin{proposition} Let $K$ and $K'$ be Lissajous knots with frequencies $(2, n_y, n_z)$ and phase shifts $(\phi_y, \phi_z)$ and   $(\phi_y', \phi_z')$ respectively. If $(\phi_y, \phi_z)$ and   $(\phi_y', \phi_z')$ are symmetric with respect to the point $(\pi/4, \pi/4)$ or the point $(\pi/4, 3 \pi/4)$ then $K$ and $K'$ are equivalent.
\label{n_x=2 symmetry} 
\end{proposition}
\begin{proof}
Suppose $(\phi_y, \phi_z)$ and   $(\phi_y', \phi_z')$ are symmetric with respect to the point $(\pi/4, \pi/4)$. Then $\phi_y'=\pi/2-\phi_y$ and $\phi_z'=\pi/2-\phi_z$. Thus
\begin{eqnarray*}
K'(-t+\pi/2)&=&(\cos(-2 t+\pi), \cos(-n_yt+n_y \pi/2+\pi/2-\phi_y), \cos(-n_z t+n_z \pi/2+\pi/2-\phi_z)\\
&=&(-\cos(2 t), (-1)^{(n_y+1)/2}\cos(n_yt+\phi_y), (-1)^{(n_z+1)/2}\cos(n_z t+\phi_z))\\
\end{eqnarray*}
which is either $K(t)$ or its mirror image $\overline{K}(t)$.

The second case follows similarly.
\end{proof}

We may now prove Theorem~\ref{first bound}.

\noindent{\bf Proof of Theorem~\ref{first bound}:}
The fundamental domain is divided into $n_x$ ``boxes'' of the form \\$[0,\frac{\pi}{n_x}]\times[k \frac{\pi}{n_x}, (k+1)\frac{\pi}{n_x}]$ for $0\le k\le n_x-1$. Within each box all the knots have the same Type II crossings and hence by Corollary~\ref{periodicity in box} there are at most $2n_y$ different knot types in that box. Since there are $n_x$ boxes we obtain at most $2 n_xn_y$ different knots.

If $n_x=2$ there is the additional rotational symmetry in each box given by Proposition~\ref{n_x=2 symmetry}. The center of each box either lies on a slanted singular line, or midway between two such lines. Moreover, one of the two boxes will be one way and the other box will be the other way. There are at most $n_y$ knot types in the box where the center of the box lies on a singular line, and there are at most $n_y+1$ knot types in the box otherwise. Thus there are at most $2 n_y+1$ knot types in total.
\hfill \qed

Our results thus far allow us to efficiently sample all Lissajous knots with a given set of frequencies $(n_x, n_y, n_z)$. We can easily pick one set of phase shifts from each region on the phase torus and Corollary~\ref{periodicity in box}, and Proposition~\ref{n_x=2 symmetry} in the case when $n_x=2$, allows us to further restrict the regions that we must sample. However, once $n_x, n_y$, and $\phi_y$ are given, the $xy$-projection has been fixed and it is natural to ask if all possible choices for $n_z$ are necessary. 
Theorem~\ref{periodicity in n_z}, which is stated in the introduction, shows  that in fact, only a finite number of values for $n_z$ are needed to produce all possible knots.

\noindent{\bf Proof of Theorem~\ref{periodicity in n_z}:}
Suppose that $K\in \L(n_x,n_y,n_z)$, $K' \in  \L(n_x,n_y,n_z+2n_xn_y)$ and that both knots have the same phase shifts. We will show first that each knot has its Type II crossings arranged the same way.

Let $(t_1, t_2)$ be a Type II crossing with parameters $(k,j)$ and let 
\begin{eqnarray*}
\Delta_{II}(n_x, n_y, n_z,\phi_z, k,j)&=&\cos(n_zt_1+\phi_z)-\cos(n_z t_1+\phi_z)\\
&=&2\sin\left(n_z \left(\frac{t_1+t_2}{2}\right)+\phi_z\right)\sin\left(n_z \left(\frac{t_1-t_2}{2}\right)\right)\\
&=&-2\sin\left(n_z\frac{j
\pi}{n_x}+\phi_z\right)
\sin\left(n_z\frac{k\pi}{n_y}\right)
\end{eqnarray*}
be the height difference between the two points on the knot directly above the crossing.

It is easy to verify that 
$$\Delta_{II}(n_x, n_y, n_z, \phi_z, k, j)=\Delta_{II}(n_z+2 n_xn_y,  n_y, n_z, \phi_z, k, j)\quad\mbox{for all}\quad k, j.$$
Thus if $n_z$ is increased by $2 n_x n_y$, not only do all Type II crossings remain unchanged, they each maintain the same height difference between upper and lower strand.

We now shift our focus to Type I crossings. Let $K$ have phase shifts $(\phi_y, \frac{n_z}{n_y}\phi_y-\varepsilon)$ and choose $\varepsilon$ small enough so that $K$ corresponds to the region just below the singular line $\phi_z=\frac{n_z}{n_y}\phi_y$. Let $K'$ correspond to the ``same'' region, that is, let $K'$ have phase shifts 
$(\phi_y, \frac{n_z+2 n_xn_y}{n_y}\phi_y-\varepsilon)$. As before, let  $(t_1, t_2)$ be a Type I crossing with parameters $(k,j)$ and let 
\begin{eqnarray*}
\Delta_{I}(n_x, n_y, n_z,\phi_y, \phi_z, k,j)&=&\cos(n_zt_1+\phi_z)-\cos(n_z t_1+\phi_z)\\
&=&2\sin\left(n_z \left(\frac{t_1+t_2}{2}\right)+\phi_z\right)\sin\left(n_z \left(\frac{t_1-t_2}{2}\right)\right)\\
&=&-2\sin\left(n_z\frac{j
\pi}{n_y}-\frac{n_z \phi_y}{n_y}+\phi_z\right)
\sin\left(n_z\frac{k\pi}{n_x}\right)
\end{eqnarray*}
be the height difference between the two points on the knot directly above the crossing.
It is easy to check that
$$\Delta(n_x, n_y, n_z, \frac{n_z}{n_y}\phi_y-\varepsilon, k, j)=
\Delta(n_x, n_y, n_z+2 n_xn_y, \frac{n_z+2 n_xn_y}{n_y}\phi_y-\varepsilon, k, j).$$
Thus $K$ and $K'$ are the same knot since both the Type I and Type II crossings are arranged the same way in each. If the phase shifts for $K$ are now changed by moving into an adjacent region, and if the phase shifts for $K'$ are changed in the same way, then the same set of crossings is changed for both $K$ and $K'$ and hence $K$ and $K'$ remain the same knot. Therefore
\begin{equation}
\label{inclusion}
\L(n_x, n_y, n_z)\subseteq \L(n_x, n_y, n_z+2 n_xn_y).
\end{equation}

According to Corollary~\ref{periodicity in box}, the pattern of knot types in each square $[0,\pi/n_x]\times[k\pi/n_x, (k+1)\pi/n_x]$, as we move from the upper left corner to the lower right corner, is periodic with period $2 n_y$. Thus if each box contains at least $2 n_y$ regions the inclusion in (\ref{inclusion}) is equality. Now the distance between successive singular lines of the type given in Equation~\ref{Type I singular lines} is 
$\frac{\pi}{\sqrt{n_y^2+n_z^2}}$ and the    distance between lines of slope 
$\frac{n_z}{n_y}$ containing opposite corners of the square  is
$\frac{(n_2+n_3)\pi}{n_1\sqrt{n_y^2+n_z^2}}$. Thus there are at least $\floor{\frac{n_y+n_z}{n_x}}$ regions in each square. Hence the inclusion in (\ref{inclusion}) is equality if   
$2 n_y\le\floor{\frac{n_y+n_z}{n_x}}.$ It is easy to check that this is true if $n_z\ge2n_xn_y-n_y$.
\hfill \qed

\begin{definition} For relatively prime integers $n_x$ and $n_y$ let $\L(n_x, n_y)=\bigcup_{n_z\in \mathbb{N}}\L(n_x,n_y,n_z)$.
\end{definition}
\begin{theorem}\label{second bound}
Let $n_x,n_y$ be relatively prime integers.  Then
$$\left|\L(n_x,n_y)\right|\leq 4n_xn_y(n_x-1)(n_y-1).$$
If furthermore $n_x=2$, then
$$\left|\L(2,n_y)\right|\leq 2(n_y-1)(2 n_y+1).$$
\end{theorem}
\begin{proof}
For fixed $n_x,n_y$ we need only consider  $2n_xn_y$ consecutive values of
$n_z$. To count the number of values that are relatively prime to both $n_x$ and $n_y$ we first subtract $2n_y$ multiples of $n_x$ that lie in
that range as well as $2n_x$ multiples of $n_y$ and then add back in the
2 multiples of $n_xn_y$.  Thus the number of possible values of
$n_z$ is bounded above by $2n_xn_y-2n_x-2n_y+2=2(n_x-1)(n_y-1)$ and
applying Theorem \ref{first bound} yields the result.
\end{proof}

\section{The Phase Torus---Fourier-$(1,1,2)$ Knots}
\label{fourier phase torus}

The phase torus of a Fourier-$(1,1,k)$ is, in general, $k+2$ dimensional although we may set any one phase shift equal to zero and drop to a $k+1$ dimension space. If $k=2$, $\phi_x=0$ and we fix $\phi_y$, then we may again think of the 2-dimensional phase torus associated to the pair $(\phi_{z,1}, \phi_{z,2})$. The singular curves are now much more complicated than in the Lissajous case, but can still be carefully described.

Suppose $K$ is a Fourier-$(1,1,2)$ knot with parameterization
\begin{eqnarray}
     x(t) & = & \cos(n_x t)\nonumber \\
     y(t) & = & \cos(n_yt+\phi_y)\label{fourier112 parameterization} \\
     z(t) & = & \cos(n_{z,1}t+\phi_{z,1})+A_{z,2}\cos(n_{z,2}t+\phi_{z,2}).\nonumber 
\end{eqnarray}

Note that by rescaling we may assume that three of the four amplitudes are $1$.

In the Lissajous case, we require that the three frequencies be pairwise relatively prime. The same proof (see \cite{BHJ1994}) can be used now to conclude that the three integers $n_x, n_y$ and $\gcd(n_{z,1},n_{z,2})$ must be pairwise relatively prime. This rules out several of the 16 cases that arise by considering all possible parities for the frequencies. Some of the remaining cases still give rise to highly symmetric knots, such as when all the frequencies are odd. In this case the knot is strongly plus amphicheiral just as in the Lissajous setting. But some of the parity cases produce knots with no apparent symmetry, suggesting that the set of Fourier-$(1,1,2)$ is much richer than the set of Lissajous knots.

We will not undertake an exhaustive analysis of the phase torus of Fourier-$(1,1,2)$ knots. Instead we offer a glimpse of the situation in the following Proposition which could be stated much more precisely. In particular, the constants in the statement of the proposition all depend on the pair of indices $(k,j)$ associated to either a Type I or II crossing. The interested reader can easily determine the constants by going through the details of the proof. Results analogous to Propositions~\ref{crossing diagonal} and \ref{n_x=2 symmetry} seem much harder.
\pagebreak
\begin{proposition}
Let  $K$ be a Fourier-$(1,1,2)$ knot with parameterization as given in \ref{fourier112 parameterization}. Then the singular curves on the phase torus are of four possible types:
\begin{enumerate}
\item Lines of the form $\phi_{z,2}=c$,
\item Lines of the form $\phi_{z,1}=c$,
\item Lines of the form $\phi_{z,2}=\pm \phi_{z,1}+c$,
\item Curves with the shape of $\sin(\phi_{z,2})=c \sin(\phi_{z,1})$
\end{enumerate}
where $c$ is a constant that, in the last case, is neither 0 nor $\pm 1$.
\end{proposition}
 \begin{proof}
  Suppose that $(t_1, t_2)$ are a pair of times that produce a double point in the $xy$-projection of $K$. Using the identity 
  $\cos x-\cos y=-2 \sin(\frac{x+y}{2})\sin(\frac{x-y}{2})$ we obtain
  $$z(t_1)-z(t_2)=-2 \sin (n_{z,1}\frac{t_1+t_2}{2}+\phi_{z,1} )\sin ( n_{z,1}\frac{t_1-t_2}{2} )
  -2 A \sin\left (n_{z,2}\frac{t_1+t_2}{2}+\phi_{z,2}\right )\sin \left ( n_{z,2}\frac{t_1-t_2}{2}\right ).$$
  We are interested in those values of $\phi_{z,1}$ and $\phi_{z,2}$ that make this difference zero.
  
  Suppose now that $(t_1, t_2)$ define a Type II crossing with indices $(k,j)$. Then
  \begin{eqnarray*}
  \frac{t_1+t_2}{2}&=&\frac{j \pi}{n_x}\\
  \frac{t_1-t_2}{2}&=&-\frac{k \pi}{n_y}
  \end{eqnarray*}
  and the crossing is singular if
  \begin{equation}
  \label{singular equation}
  \sin \left (\frac{n_{z,1} j \pi}{n_x}+\phi_{z,1}\right )\sin \left(\frac{n_{z,1} k \pi}{n_y}\right )=
 -A\sin\left(\frac{n_{z,2} j \pi}{n_x}+\phi_{z,2}\right )\sin \left(\frac{n_{z,2} k \pi}{n_y}\right ).
 \end{equation}
  We are now led to several cases.
  
  \noindent{\bf Case I: $n_y |   n_{z,1 } k$}\\
   If $n_y$ divides $ n_{z,1}$ then we must have that $\sin\left(\frac{n_{z,2} j \pi}{n_x}+\phi_{z,2}\right )=0$ since $k<n_y$  and $n_y, n_{z,1}$ and $n_{z,2}$ cannot have a common factor. This means that
  $$\phi_{z,2}=m \pi-\frac{n_{z,2} j \pi}{n_x}$$
  for some integer $m$.
  
    \noindent{\bf Case II: $n_y | n_{z,2} k$}\\
     This is similar to Case I leading to
     $$\phi_{z,1}=m \pi-\frac{n_{z,1} j \pi}{n_x}$$
  for some integer $m$.
  
  If the first two cases do not occur, then we may rewrite Equation~\ref{singular equation} as 
  \begin{equation*}
  \sin \left (\frac{n_{z,1} j \pi}{n_x}+\phi_{z,1}\right )=
 C\sin\left(\frac{n_{z,2} j \pi}{n_x}+\phi_{z,2}\right )
 \end{equation*}
 where
 $$C= -A   \sin \left(\frac{n_{z,2} k \pi}{n_y}\right )/\sin \left(\frac{n_{z,1} k \pi}{n_y}\right ).$$
 
 \noindent{\bf Case III: $|C|=1$}\\
  In this case we must have 
 \begin{equation}
 \left (\frac{n_{z,1} j \pi}{n_x}+\phi_{z,1}\right ) \pm\left (\frac{n_{z,2} j \pi}{n_x}+\phi_{z,2}\right )= m \pi
 \label{sum or difference}
 \end{equation}
  for some integer $m$, where the parity of $m$ depends on the sign of $C$ and whether we are forming a sum or difference in Equation~\ref{sum or difference}. Thus $\phi_{z,2}=\pm \phi_{z,1}+c$ for some constant $c$.
 
  \noindent{\bf Case IV: $|C|\ne1$}\\
   In this case we are left with a translate of the curve 
  $$\sin(\phi_{z,1})=C \sin(\phi_{z,2}).$$
  This is an interesting curve which, at first glance,  appears much like a sine curve. It is oriented either vertically or horizontally depending on the value of $|C|$.
  
  The analysis of a Type I crossing is similar and is left to the reader.
 \end{proof}

\begin{figure}
    \begin{center}
    \leavevmode
    \scalebox{.50}{\includegraphics{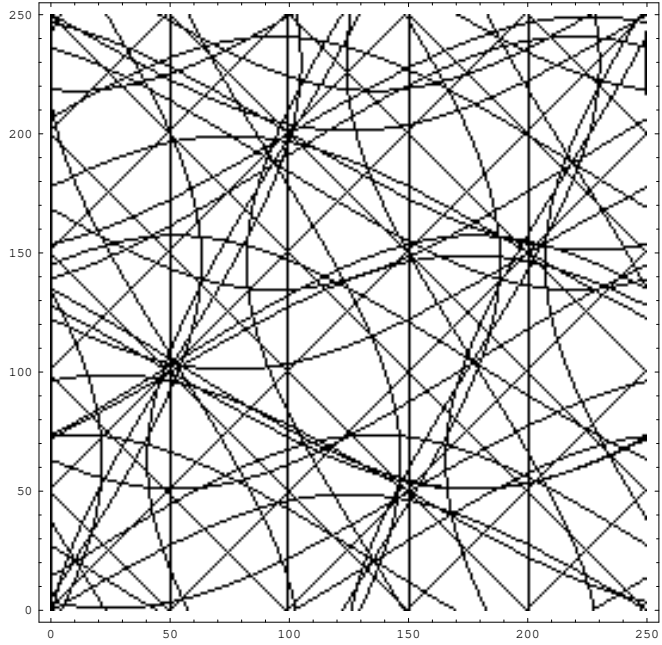}}
    \end{center}
\caption{The phase torus for the Fourier-$(1,1,2)$ knot with $n_x=5, n_y=6, n_{z,1}=1, n_{z, 2}=2, \phi_x=0, \phi_y=\pi/4$ and $A_{z,1}=1$,   shown for $0\le \phi_{z,1}\le \pi$ and $ 0\le \phi_{z,2}\le \pi$.}
\label{5612phaseTorus}
\end{figure}

In Figure~\ref{5612phaseTorus} we give an example showing a 250 pixel by 250 pixel bitmap image of the phase torus for a specific set of parameters.  Even with relatively small frequencies, one can begin to appreciate the difficulty of systematically sampling each region of the phase torus for an arbitrary Fourier-$(1,1,2)$ knot.

\section{2-Bridge Knots}
\label{2-bridge knots}

Every 2-bridge knot can be classified by a pair of relatively prime integers $(p,q)$ such that $p$ is odd and $0<q<p$. We will often write the pair $(p,q)$ as the fraction $p/q$. If $K_{p/q}$ and $K_{p'/q'}$ are two 2-bridge knots with corresponding fractions $p/q$ and $p'/q'$ then they are equivalent knots if and only if $p=p'$ and 
$\pm q'q^{\pm 1} \equiv 1$ mod $p$. 
The reader is referred to \cite{BZ:2003} for details.

If $K$ is a Fourier~$(1,1,2)$ knot with $n_x=2$ then $K$ is a 2-bridge knot. We may recover the fraction $p/q$ from the Lissajous projection in the $xy$-plane as follows. This projection is a 4-plat diagram. As we move in the $x$-direction from left to right we see a single Type I crossings on the $x$-axis, then a pair of Type II crossings which are symmetric with respect to the $x$-axis, then another Type I crossing on the $x$-axis, and so on.  Let $\eta_1, \eta_2, \dots$ be the signs of the Type I crossings from left to right along the $x$-axis. Let $\{ \varepsilon_1^1,  \varepsilon_1^2\}, \{ \varepsilon_2^1,  \varepsilon_2^2\}, \dots$ be the signs of the pairs of Type II crossings  from left to right. Proceeding in a fashion similar to that given on pages 300--303 in \cite{Rolfsen}, we obtain that  $p/q$ is given by the continued fraction
\begin{equation}\label{continued fraction}
p/q=[\eta_1,\varepsilon_1^1+\varepsilon_1^2,\eta_2,\varepsilon_2^1+\varepsilon_2^2,\ldots,\eta_{n_y}]=\eta_1+\cfrac{1}{\varepsilon_1^1+\varepsilon_1^2 +
\cfrac{1}{\eta_2+\dots +\cfrac{1}{\eta_{n_y}}
}}
\end{equation}
Note that  if $K$ is Lissajous then it  is rotationally symmetric with respect to the $x$-axis and each pair of Type II crossings has the same sign. In this case each $\varepsilon_i^1+\varepsilon_i^2$ can be replaced with $2\varepsilon_i^1$.
Using this formula,  it is easy to determine the 2-bridge knot given by a Fourier-$(1,1,2)$ representation with $n_x=2$. Hence, when we sample Lissajous and Fourier-$(1,1,2)$ knots with $n_x=2$, even if we obtain knots with hundreds of crossings, it is a simple matter to distinguish them.

Since every Lissajous knot with $n_x=2$ is 2-bridge, a good question is: What 2-bridge knots are Lissajous with $n_x=2$? As mentioned in the Introduction, every Lissajous knot is either strongly plus amphicheiral, or 2-periodic and linking its axis of rotation once. It is known that a 2-bridge knot cannot be strongly plus amphicheiral \cite{HK:1979}. It is also known (and will be shown below) that every 2-bridge knot is 2-periodic, but may or may not link its axis of rotation once. The following theorem makes it easy to identify which 2-bridge knots {\sl might} be Lissajous.

\begin{theorem} Let $K$ be a 2-bridge knot. Then the following are equivalent.

1) $K$ has a  symmetry of period 2 with axis $A$ such that $A$ is disjoint from $K$ and $|lk(A,K)|=1$

2) $\Delta_K (t)$ is a square mod 2

3) $\Delta_K(t)\equiv$ 1 mod 2.

\end{theorem}
\begin{proof} As already mentioned in the introduction, it follows from a result of Murasugi  (see \cite{M:1971}) that 1) implies 2) and clearly 3) implies 2). We must show that 2) implies both 1) and 3).

Suppose $K$ is a 2-bridge knot given by the pair of relatively prime integers $(p,q)$ with $p$ odd and $0<q<p$. There is a unique continued fraction expansion 
$$p/q=[2a_1, -2a_2, \dots, (-1)^{n+1}2a_n]$$
where each of the partial quotients $(-1)^{i+1}2a_i$ is even. Corresponding to this expansion is a Seifert surface made from plumbing together twisted bands as shown in Figure~\ref{plumbed bands}. Notice that the Seifert surface, and hence $K$ is rotationally symmetric   around the axis $A$. Thus every 2-bridge knot has a  symmetry of period 2 with axis disjoint from the knot. However, the linking number of $A$ and $K$ need not be $\pm1$ in general. From the plumbing picture we also see that the number of bands, $n$,  must be even in order to get a knot. If $n$ is odd we obtain a 2-bridge link (of two components).
\begin{figure}
\psfrag{a}{$a_1$}
\psfrag{b}{$a_2$}
\psfrag{S}{$S$}
\psfrag{A}{$A$}
\psfrag{c}{$a_n$}
    \begin{center}
    \leavevmode
    \scalebox{.85}{\includegraphics{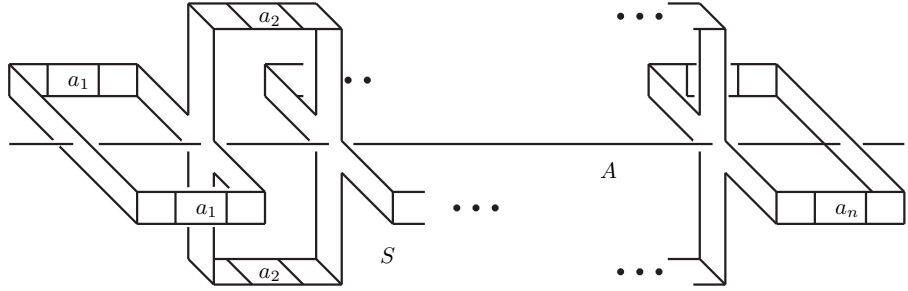}}
    \end{center}
\caption{The Seifert surface $S$ and axis $A$ for the 2-bridge knot $K_{p/q}$. Each $a_i$ represents $a_i$ right handed half-twists in the band.}
\label{plumbed bands}
\end{figure}

The axis $A$ meets the Seifert surface $S$ transversely in $n+1$ points and we may compute the linking number of $A$ and $K$ by counting the signed intersection points. Let $\epsilon_i$ be the sign of the intersection point that occurs between band $i$ and $i+1$. Let $\epsilon_0$ be the sign of the left-most intersection point in the figure and choose orientations so that $\epsilon_0=1$. It is easy to see that
$\epsilon_{i+1}=\epsilon_i$ if $a_i$ is odd and $\epsilon_{i+1}=-\epsilon_i$ if $a_i$ is even. Thus the sequence $\{a_1, a_2, \dots, a_n\}$ determines the sequence $\{\epsilon_0, \epsilon_1, \dots, \epsilon_n\}$ which in turn determines the linking number between $A$ and $K$.

From the Seifert surface we may obtain the Seifert matrix $V$ and compute the Conway polynomial $\nabla(z)=\det (t^{-1/2}V-t^{1/2}V^T)$. It is a straightforward calculation to show that
$$\nabla_K(z)=\begin{array}{c}(\begin{array}{cc}1&0\end{array}) \\ \\\end{array} \left( \begin{array}{cc} -a_1 z&1\\1&0\end{array}\right )
\left ( \begin{array}{cc} -a_2 z&1\\1&0\end{array}\right )
\dots
\left ( \begin{array}{cc} -a_n z&1\\1&0\end{array}\right )\left (\begin{array}{c}1\\0\end{array} \right)$$
(See page 207 of \cite{C:2004}.)

If $n=2$ the Conway polynomial is $\nabla(z)=1+a_1a_2 z^2$ and the Alexander polynomial $\Delta(t)=\nabla(t^{1/2}-t^{1/2})=a_1a_2 t^{-1}+(1-2a_1a_2)+a_1a_2t$. In general, an Alexander polynomial of the form 
$$\Delta(t)=b_0+b_1(t+t^{-1})+b_2(t^2+t^{-2})+\dots+b_m(t^m+t^{-m})$$ is a square mod 2 if and only if $b_{2k+1}\equiv 0$ mod 2 for $k=0, 1, \dots$. Thus 
$a_1a_2 t^{-1}+(1-2a_1a_2)+a_1a_2t$ is a square mod 2 if and only if at least one of $a_1$ or $a_2$ is even. But if this is the case it follows that $\sum \epsilon_i=\pm 1$. It is also true that $a_1a_2 t^{-1}+(1-2a_1a_2)+a_1a_2t$ is a square mod 2 if and only if it is equal to 1 mod 2.

Suppose now that $n>2$ and that $K$ is a knot with $\Delta_K(t)$ a square mod 2.  The first thing to show is that $a_i$ is even for at least one value of $i$. If $a_i$ were odd for every $i$, then replace each $a_i$ with $-1$.  This does not change any $a_i$ mod 2, and hence  does not change $\nabla_K(z)$ mod 2 or $\Delta_K(t)$ mod 2.  But if $a_i=-1$ for all $i$, we can prove by induction on $n$ that 
$$\Delta(t)=1-t+t^2-t^3+\dots+t^n$$
and it follows that $\Delta(t)$ is not a square mod 2. Thus at least one $a_i$ is even. Replacing this $a_i$ with zero transforms $K$ into a knot $J$ with the same Alexander polynomial mod 2 but with a Seifert surface having two fewer bands. Proceeding by induction on the number of bands, we have that lk$(J,A)=\pm 1$. But now because $a_i$ is even, $K$ must also link $A$ once. If we begin instead with the assumption that  $\Delta_K(t)\equiv 1$ mod 2, the same argument will work.
\end{proof}

\section{Sampling Lissajous and Fourier Knots}
\label{sampling}
Using the results of Sections~\ref{lissajous phase torus} we are now in a position to efficiently sample Lissajous knots. In the case where $n_x=2$ we obtain 2-bridge knots and can take advantage of this to compare knots in our sample. For the more general case of Fourier knots, we have not carried out a complete analysis of the phase torus, a task that seems much more difficult. Hence, we have not attempted to rigorously sample Fourier knots, but instead have relied on two methods, either random sampling or an algorithm which  first ``draws''  a bitmap image of the phase torus (as in Figure~\ref{5612phaseTorus}) and then  picks one point from each ``white'' region. This latter approach is fraught with difficulty since, for example,  some white regions may be smaller than a single pixel and be missed. Our samples naturally fall into four cases which we describe in turn in this section.

\subsection{Lissajous Knots with $2=n_x<n_y<n_z$}
\label{sampling Lissajous n_x=2}
We have determined all knots in $\L(2, n_y)$ for $3\le n_y\le 105$. For a given value of $n_y$ we let $n_z$  run from $3 n_y+2$ to $7n_y$. These values of $n_z$ are sufficient to guarantee that we obtain all possible knots in  $\L(2, n_y)$. Since each of these knots is 2-bridge we were able to use  Equation~\ref{continued fraction}  to identify the associated pair $(p,q)$ and thus compare knots in the output.  The total number of knots in $\L(2, n_y)$ is given in 
Table~\ref{size of L(2,n_y)} for each value of $n_y$. It is interesting to compare these numbers with the upper bound given by Theorem~\ref{second bound}. Depending on $n_y$, the actual number of knots found  is roughly between 5 and 10 per cent of the upper bound. The discrepancy is almost certainly due to the presence of huge numbers of unknots.  The $xy$-projection of a Lissajous knot with $n_x=2$ and $n_y=99$ has $(2)(2)(99)-2-99=295$ crossings, and knots in $\L(2,99)$ have crossing numbers ranging from 5 to 293. Of course the bound of 78008 given by Theorem~\ref{second bound} for $n_x=2$ and $n_y=99$ is well below the upper bound of $2^{295}$ obtained by considering all possible crossing arrangements!

\begin{table}
\begin{center}
\begin{tabular}{|cc|cc|cc|cc|}
\hline
$n_y$ & $|\L(2, n_y)|$ &$n_y$ & $|\L(2, n_y)|$ & $n_y$ & $|\L(2, n_y)|$ &  $n_y$ & $|\L(2, n_y)|$ \\
\hline
3&3&29&645&55&1854&81&3761\\
5&11&31&737&57&1727&83&5805\\
7&28&33&533&59&2859&85&4654\\
9&37&35&684&61&3062&87&4195\\
11&78&37&1075&63&1946&89&6707\\
13&109&39&772&65&2639&91&5647\\
15&93&41&1339&67&3708&93&4805\\
17&203&43&1473&69&2593&95&5892\\
19&258&45&904&71&4191&97&7984\\
21&195&47&1782&73&4433&99&5208\\
23&390&49&1688&75&2584&101&8699\\
25&390&51&1365&77&3933&103&9036\\
27&387&53&2287&79&5248&105&4425\\
\hline
\end{tabular}

\caption{The number of distinct Lissajous knots with $n_x=2$ as a function of $n_y$.}
\label{size of L(2,n_y)}
\end{center}
\end{table}

The total number of knots in Table~\ref{size of L(2,n_y)} is 135061, far too many to describe one by one. However, in Tables~\ref{p/q for n_y up to 15 part 1}--\ref{p/q for n_y up to 15 part 4} we list all knots in $\L(2, n_y)$, grouped by crossing number, for $3\le n_y\le 15$. Several interesting things can be seen in these tables. The same knot often appears in $\L(2, n_y)$ for many different values of $n_y$. For example $K_{7/2}$ (which is the twist knot $5_2$ in \cite{Rolfsen}) appears in every column of Table~\ref{p/q for n_y up to 15 part 1}. In fact, $K_{7/2}\in \L(2, n_y)$ for  $3\le n_y\le 105$.  This is also true for $K_{9/2}$. The knot $K_{15/4}$ first appears for $n_y=3$, misses a few values of $n_y$, and then is contained in $\L(2, n_y)$ for $23\le n_y\le 105$.  Similar patterns hold for the other small-crossing knots suggesting that if $K\in \L(2, n_y)$ for some $n_y$ then there exists $N$ such that $K\in \L(2, n_y)$ for all $n_y\ge N$. A second observation is that several small crossings knots are already conspicuously absent. In particular, there are exactly four 8-crossings knots with Alexander polynomial congruent to 1 mod 2 (and hence possibly Lissajous). These are $K_{17/4}, K_{23/7}, K_{25/9}$ and $K_{31/12}$, only one of which, $K_{31/12}$, appears to be Lissajous. While 
Tables~\ref{p/q for n_y up to 15 part 1}--\ref{p/q for n_y up to 15 part 4} display only a small fraction of our total sample, it is in fact true that the other three 8-crossing knots  do not appear for any $n_y$ up to 105.

\begin{question} Does there exist a 2-bridge knot $K$ with $\Delta_K(t)\equiv 1$ mod 2 that is not Lissajous (with or without one frequency equal to 2)? In particular, are any of the 8-crossing  2-bridge knots  $K_{17/4}, K_{23/7}$ or $ K_{25/9}$ Lissajous?
\end{question}

In Table~\ref{2-bridge numbers} we list the numbers of  2-bridge knots, 2-bridge knots with Alexander polynomial congruent to 1 mod 2, and finally, the number of these that are Lissajous knots with $n_x=2$ and $3\le n_y\le 105$. The table has entries for each crossing number from 3 to 16. Very quickly we see that many 2-bridge knots with the required symmetry are not Lissajous, at least not with $n_x=2$ and $3\le n_y\le 105$. It seems unlikely that choosing $n_y>105$ will yield more 2-bridge knots in the 3--16 crossing range. On the other hand, perhaps letting the even frequency be more than 2 will yield more 2-bridge knots with small crossing number. We examine this further in Section~\ref{sampling Lissajous n_x>2}.

\begin{table}[htdp]
\begin{center}
\begin{tabular}{|c|cccccccccccccc|}
\multicolumn{15}{c}{crossing number}\\
\hline
&3&4&5&6&7&8&9&10&11&12&13&14&15&16\\
\hline
2-bridge&1&1&2&3&7&12&24&45&91&176&352&693&1387&2752\\
$\Delta (t) \equiv 1$&0&0&1&1&2&4&8&13&26&51&97&185&365&705\\
$\L(2,n_y)$&0&0&1& 1& 2 & 1 & 3 & 4 & 8 & 5 & 9 & 7 & 15 & 15\\
\hline

\end{tabular}
\end{center}
\caption{The number of  2-bridge knots, 2-bridge knots with Alexander polynomial congruent to 1 mod 2, and the number of these that are Lissajous with $n_x=2$ and $3\le n_y\le 105$, as a function of crossing number.}
\label{2-bridge numbers}

\end{table}%

In Tables~\ref{2-bridge Lissajous knots 3 to 16 part 1}--\ref{2-bridge Lissajous knots 3 to 16 part 2} we list all 2-bridge knots with crossings from 3 to 16 which are Lissajous knots with $n_x=2$ and $3\le n_y \le 105$.  Here the knot name, as defined in \cite{HTW:1998} and used in {\it Knotscape} \cite{HT:1998} appears in the first column. Following that we give the 2-bridge defining fraction $p/q$  and then the frequency and phase shift data. For each knot, the given value of $n_y$ is minimal. However, since our search let $n_z$ run from $3 n_y+2$ to $7n_y$, it might be possible for a given knot to be represented with a smaller value of $n_z$.

As a check against errors, we took all the 2-bridge knots in the data set that have Lissajous diagrams with less than 50 crossings (the built-in limit for {\it Knotscape}) and crossing number less than 17, and looked them up in the {\it  Knotscape} table of knots in two different ways. First we  converted their Lissajous diagrams to Dowker-Thistlethwaite code (the input format for {\it Knotscape}) and then used the ``Locate in Table'' feature. Next we converted the defining fraction $p/q$ into  DT code and again used the ``Locate in Table'' routine. Happily, the results matched.

\subsection{Lissajous Knots with $2<n_x<n_y<n_z$}
\label{sampling Lissajous n_x>2}

Our goal in this section is to simply find as many Lissajous knots in the 3--16 crossing range as we can. We may still use the results of Section~\ref{lissajous phase torus} to efficiently sample Lissajous knots with all frequencies greater than 2, but it is more difficult to tabulate the output. This is because even with relatively small frequencies, very large crossing number knots can result, and we can no longer use the classification of 2-bridge knots to sort them out. Therefore, we limited ourselves to producing diagrams with at most 49 crossings, the limit of what can be input to {\it Knotscape}.  Assuming that $2<n_x<n_y$, and that $\gcd(n_x, n_y)=1$, we are left with the following $(n_x, n_y)$ pairs:
$$\{(3,4), (3,5),(3,7),(3,8),(3,10), (4,5),(4,7),(5,6)\}.$$ For each of these pairs we let $n_z$ run from $2n_xn_y-n_x-n_y$ to $4n_xn_y-n_x-n_y-1$, a range sufficient to produce all possible Lissajous knots. We obtained a total of 6352 knots of which Knotscape identified 1428 as unknots. The remaining 4924 knots fell into four categories:
\begin{enumerate}
\item knots identified as composites by {\sl Knotscape},
\item knots which {\sl Knotscape} located in the Hoste-Thistlethwaite-Weeks table,
\item knots which {\sl Knotscape} simplified to alternating projections with more than 16 crossings, and 
\item knots which {\sl Knotscape} simplified to nonalternating projections with more than 16 crossings.
\end{enumerate}

In Table~\ref{composite Lissajous knots} and Tables~\ref{Lissajous knots part 1}--\ref{Lissajous knots part 2} we list all knots in the first two categories. We note that while {\sl Knotscape} can identify a knot as a composite, it identifies the summands only up to mirror image. In order to properly identify the composites in Table~\ref{composite Lissajous knots} we compared their Jones polynomials to the Jones polynomials of all possible composites using the given summands or their mirror images in all possible ways.

\begin{table}[h!]
\begin{center}
\begin{tabular}{|l|cccccc|}
\hline
knot& $n_x$ & $n_y$ & $n_z$ & $\phi_x$ & $\phi_y$ & $\phi_z$\\
\hline
  
 $3a1\# 3a1$&  3 & 4 & 23 & 0 & 0.25210 & 1.84229 \\
 $3a1\# \overline{3a1}$&3 & 5 & 29 & 0 & 0.23099 & 2.91059 \\
  $5   a   1 \#  5   a   1$&3 & 7 & 50 & 0 & 0.50522 & 1.58916 \\
  $5   a   1 \#  \overline{5   a   1}$&3 & 5 & 29 & 0 & 0.26179 & 1.83259 \\
$6   a   1  \#  6   a   1  $&4 & 5 & 37 & 0 & 0.18699 & 2.95459 \\
  $6   a   3  \#  6   a   3$&3 & 8 & 47 & 0 & 0.23799 & 0.80919 \\
  $6   a   3  \#  \overline{6   a   3} $&3 & 5 & 29 & 0 & 0.29259 & 0.75459 \\
  $3   a   1  \#  3   a   1    \#5   a   1$&4 & 5 & 39 & 0 & 0.16064 & 2.19554 \\
  $3   a   1  \#  3   a   1    \#\overline{5   a   1}$&4 & 7 & 55 & 0 & 0.13934 & 2.21684 \\
$  3   a   1   \# 3   a   1 \#   8   a   2 $&5 & 6 & 59 & 0 & 0.11116 & 2.40211 \\
  $5   a   1  \#  5   a   1  \# \overline{ 5   a   1}$&4 & 7 & 55 & 0 & 0.15201 & 1.41878 \\
  $6   a   3  \#  6   a   3   \# \overline{6   a   3} $&4 & 7 & 55 & 0 & 0.16468 & 0.62071 \\
  $3   a   1  \#  3   a   1 \#   3   a   1 \#   3   a   1  $&5 & 6 & 59 & 0 & 0.10149 & 1.78345 \\

  \hline
\end{tabular}
\end{center}
\caption{Small-crossing composite Lissajous knots.  A bar over a knot name indicates mirror image. Knot names are as in {\sl Knotscape}.}
\label{composite Lissajous knots}
\end{table}%

The third category cannot include knots in the Hoste-Thistlethwaite-Weeks table and we make no attempt to list them here. The fourth category {\it might} have included knots with 16 or less crossings that {\sl Knotscape} simply failed to simplify correctly.  To investigate this we first computed the Jones polynomial  of each knot and eliminated  knots whose Jones polynomial had a span of 17 or more. (Recall that the crossing number of a knot is bounded below by the span of the Jones polynomial.) This left a total of 78 knots. Of these, only 5 shared the same Jones polynomial with prime knots having less than 17 crossings and furthermore having an Alexander polynomial that is a square mod 2. In each of these five cases either the Alexander polynomial or the Kauffman 2-variable polynomial was sufficient to show that the knots did indeed have crossing numbers of 17 or more.

Thus, barring clerical errors, Table~\ref{composite Lissajous knots} and Tables~\ref{Lissajous knots part 1}--\ref{Lissajous knots part 2} provide a complete list of all Lissajous knots with $x$ and $y$ frequencies of $(3,4), (3,5),(3,7),(3,8),(3,10), (4,5),(4,7)$ or $(5,6)$ which are either composite, or prime  with 16 or less crossings.

As mentioned in the introduction, there are exactly three prime knots with 12 or less crossings that are strongly plus amphicheiral: 10a103 ($10_{99})$, 10a121 ($10_{123}$), and 12a427. The knots 10a103 and 12a427 are Lissajous and are listed in Table~\ref{Lissajous knots part 1}. A natural question is,
\begin{question} Is the strongly plus amphicheiral knot 10a121 Lissajous?
\end{question}

The knot 10a121 is one member of a family of knots known as {\it Turks Head} knots. These knots are conjectured to not be Lissajous by Przytycki. See \cite{P}.

It is easy to see that every composite knot of the form $K\#\overline{K}$ is strongly plus amphicheiral while composites of the form $K\# K$ are 2-periodic and link their axis of rotation once. Several knots of this form appear in  Table~\ref{composite Lissajous knots}. Thus another good question is,
\begin{question}Is every composite knot of the form $K\#K$ or  $K\#\overline{K}$ Lissajous?
\end{question}

\subsection{Fourier-$(1,1,2)$ Knots with $2=n_x<n_y$}
\label{sampling Fourier n_x=2}
Rather than trying  to algorithmically choose one point in each region of the phase torus for a Fourier-$(1,1,2)$ knot, we chose instead to randomly sample points from the phase torus. Fixing $n_x=2$, $\phi_x=0$ and $A_{z,1}=1$, we then let $n_y$ take on odd values from $3$ to $99$. For each value of $n_y$ the remaining parameters were then chosen at random such that:

\begin{center}
\begin{tabular}{c}
$\phi_y=\frac{k}{7}\pi\mbox{, }k\in\{1,2,3,4,5,6\}$\\
 $0<n_{z,1}<n_{z,2}<301$\\
 $0\le \phi_{z,1}\le \pi$\\
 $0\le \phi_{z,2}\le 2 \pi$\\
 $0\le A_{z,2}\le 2$\\
 \end{tabular}
 \end{center}
 
 For each value of $n_y$, random sampling in batches of 10000 took place until no new knots were found. If a knot was produced that had already been found, the one with the lexicographically smallest set $\{n_x, n_y, n_{z,1}, n_{z,2}\}$ was kept. This tended to produce knots with fairly small values of $\{n_x, n_y, n_{z,1}$ but with $n_{z,2}$ often in the hundreds. Furthermore, only knots with less than 17 crossings were kept in the sample.
 
\begin{figure}[h!]
    \begin{center}
    \leavevmode
    \scalebox{.70}{\includegraphics{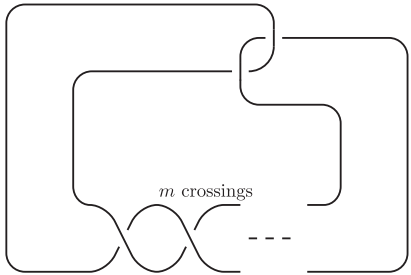}}
    \end{center}
\caption{The twist knot $T_m$.}
\label{twistKnot}
\end{figure}

 After a modest amount of searching, we turned up all 2-bridge knots with 14 or less crossings, and nearly all 15 and 16-crossing ones as well. (We found 1386 out of 1387 15-crossing knots and 2731 out of  2752 16-crossing knots.) We believe the following conjecture is reasonable.
 \begin{conjecture}
 \label{all 2-bridge are fourier (1,1,2) conjecture}
 Every 2-bridge knot can be expressed as a Fourier-$(1,1,k)$ knot with $n_x=2$ and $k\le 2$.
 \end{conjecture}
Additional evidence for this conjecture is provided by the twist knots. The twist knot $T_m$, which is the 2-bridge knot $K_\frac{2 m+1}{2}$, is shown in Figure~\ref{twistKnot}. The mirror image of $T_m$  is the twist knot $T_{-1-m}$. Thus it suffices to consider $m>1$.  It is shown in \cite{HZ2006} that $T_m$ is Lissajous if and only if $m\equiv 0$ mod 4 or $m\equiv 3$ mod 4. If this is not the case, the knot does not have the required symmetry to be Lissajous. However, in these cases, the following examples show that $K_m$ is a Fourier-$(1,1,2)$ knot. Thus all twist knots are Fourier-$(1,1,k)$ knots with $k\le 2$.

\begin{theorem} Twist knots which are not Lissajous may be expressed as Fourier knots as follows.
\begin{enumerate} 
\item The twist knot $T_{4 n+1}$ can be expressed as the Fourier-$(1,1,2)$ knot with $n_x=2,\phi_x=0,  n_y=8n+3, \phi_y=1/2, n_{z,1}=2, \phi_{z,1}=\pi/4, n_{z,2}=8n+1, \phi_{z,2}=\frac{8n+1+(8n+5)\pi}{2 (8n+3)}$ and $A_{z,2}=1$ for all $n\ge 1$. 

\item The twist knot $T_{2n}$ can be expressed as the Fourier-$(1,1,2)$ knot with $n_x=2,\phi_x=0,  n_y=2 n+1, \phi_y=1/2, n_{z,1}=2, \phi_{z,1}=\pi/4, n_{z,2}=2n+3, \phi_{z,2}=\frac{2n+3-3\pi}{2 (2n+1)}$ and $A_{z,2}=1$ for all $n\ge 1$.
\end{enumerate}
\end{theorem}
The proof is similar to the proof of Theorem~4 given in \cite{HZ2006} and relies on very carefully determining the sign of each crossing in the diagram. The details are quite long and not particularly insightful. We leave this as a rather complicated exercise for the reader.

 Our sample of all 2-bridge knots to 16 crossings expressed as Fourier-$(1,1,2)$ knots is too large to reproduce here. Instead, in Tables~\ref{2-bridge 3 to 10 as fourier(1,1,2), part 1}--\ref{2-bridge 3 to 10 as fourier(1,1,2), part 2}
 we list all 2-bridge knots to 10 crossings with associated Fourier data. To generate this table we again undertook a random sample but this time sharply reduced the range of the parameters. In particular, we kept all amplitudes equal to one, set $\phi_y=\pi/4$, and only allowed $z$-frequencies as large as 10. An interesting variation on Conjecture~\ref{all 2-bridge are fourier (1,1,2) conjecture} would be to require that all amplitudes are 1. Knots appearing in Tables~\ref{2-bridge 3 to 10 as fourier(1,1,2), part 1}--\ref{2-bridge 3 to 10 as fourier(1,1,2), part 2} which are known to be Lissajous are shown in boldface, while those that have Alexander polynomials congruent to 1 mod 2, and hence {\it might} be Lissajous, are shown in italics.

\subsection{Fourier-$(1,1,2)$ Knots with $2<n_x<n_y$}
\label{sampling Fourier n_x>2}

We made only a modest attempt to sample Fourier-$(1,1,2)$ knots with $x$ and $y$ frequencies greater than two. Rather than sampling at random as in Section~\ref{sampling Fourier n_x=2}, we now chose one sampling point from each region of the phase torus by first creating a bitmap image as in Figure~\ref{5612phaseTorus} and then taking the centroid of each white region. Sometimes the centroid fell outside of the region and in this case an arbitrary point of the region was selected. Because of the large crossing numbers that result, and the consequent difficulty in identifying these knots, we again restricted our sample to  $x$ and $y$ frequencies of $(3,4), (3,5),(3,7),(3,8),(3,10), (4,5),(4,7)$ or $(5,6)$. We further restricted the $z$-frequencies to be less than 15 and somewhat arbitrarily fixed all amplitudes at 1. Using {\sl Knotscape} to identify the resulting knots, and keeping only knots with 16 or less crossings, we found several thousand prime knots. In Tables~\ref{non 2-bridge 3 to 10 as fourier(1,1,2), part 1}--\ref{non 2-bridge 3 to 10 as fourier(1,1,2), part 3} we list all of these with 10 or less crossings  which are not 2-bridge and thus listed in Tables~\ref{2-bridge 3 to 10 as fourier(1,1,2), part 1}--\ref{2-bridge 3 to 10 as fourier(1,1,2), part 2}. All knots through 9 crossings were found, and all but 20 alternating 10-crossings knots were found. We suspect that limiting the $z$-frequencies to less than 15 is a severe restriction.

We did however, find all torus knots up to 16 crossings. It is shown by Kauffman in \cite{K} that every torus knot is a Fourier-$(1,3,3)$ knot. Interestingly, we found that, up to 16 crossings,  the torus knot $T_{p,q}$ can be represented by a Fourier-$(1,1,2)$ knot with $n_x=p$ and $n_y=q$. Table~\ref{torus knots} lists these results. The data suggests the following conjectures (which we have verified for a large number of $(p,q)$ pairs, and hope to prove in a later paper). 
\begin{conjecture} The torus knot $T_{2,q}$ can be represented as a Fourier-$(1,1,2)$ knot with frequencies $n_x=2, n_y=q, n_{z,1}=2$ and $n_{z,2}=q-2$ and phase shifts $\phi_x=0, \phi_y=\pi/4, \phi_{z,1}=\pi/2$ and $\phi_{z,2}=\pi/4$.
\end{conjecture}

\begin{conjecture} The torus knot $T_{p,q}$, with $0<p<q$,  can be represented as a Fourier-$(1,1,2)$ knot with frequencies $n_x=p, n_y=q, n_{z,1}=p$ and $n_{z,2}=q-p$.
\end{conjecture}

It would be interesting to undertake a large-scale sampling of Fourier-$(1,1,2)$ knots with $x$ and $y$ frequencies greater than two to see if {\it every} knot with 16 or less crossings turns up. Such a study might shed light on the  following question,

\begin{question}Is there a knot which cannot be expressed as a Fourier-$(1,1,k)$ knot for $k\le 2$?
\end{question}

\pagebreak
\section{Tables of Lissajous and Fourier Knots}
\begin{table}[h!]
\begin{center}
\begin{tabular}{|l|c|c|c|c|c|c|c|}
\hline
&$\L(2,3)$&$\L(2,5)$&$\L(2,7)$&$\L(2,9)$&$\L(2,11)$&$\L(2,13)$&$\L(2,15)$\\
\hline
cr&$p/q$&$p/q$&$p/q$&$p/q$&$p/q$&$p/q$&$p/q$\\
\hline
5 & 7/2 & 7/2 & 7/2 & 7/2 & 7/2 & 7/2 & 7/2 \\
\hline
6 & 9/2 & 9/2 & 9/2 & 9/2 & 9/2 & 9/2 & 9/2 \\
\hline
7 & 15/4 & 17/5 & 15/4 & 17/5 & 15/4 & 15/4 & 17/5 \\
7 &   &   & 17/5 &   & 17/5 & 17/5 &   \\
\hline
8 &   &   &   &   & 31/12 &   &   \\
\hline
9 &   &   & 15/2 & 31/7 & 15/2 & 15/2 & 15/2 \\
9 &   &   & 31/7 &   & 31/7 & 31/7 & 31/7 \\
9 &   &   &   &   & 31/11 &   &   \\
\hline
10 &   & 17/2 & 17/2 & 17/2 & 55/12 & 17/2 & 55/12 \\
10 &   & 49/20 & 49/20 & 57/16 & 57/16 & 49/20 & 57/16 \\
10 &   & 57/16 & 55/12 &   &   & 57/16 &   \\
10 &   &   & 57/16 &   &   &   &   \\
\hline
11 &   & 65/14 & 73/16 & 49/9 & 49/9 & 41/13 & 49/9 \\
11 &   & 73/16 &   & 65/14 & 73/16 & 49/9 & 73/16 \\
11 &   & 97/26 &   & 73/16 &   & 71/20 &   \\
11 &   &   &   &   &   & 73/16 &   \\
11 &   &   &   &   &   & 97/26 &   \\
\hline
12 &   & 121/32 & 169/50 & 167/46 & 121/32 & 167/46 & 169/50 \\
12 &   &   &   &   & 169/50 &   &   \\
\hline
13 &   & 209/56 & 239/71 &   & 23/2 & 71/11 & 23/2 \\
13 &   &   &   &   & 71/11 & 209/56 & 71/11 \\
\hline
14 &   &   & 25/2 & 25/2 & 407/119 & 25/2 & 407/119 \\
14 &   &   & 89/36 & 89/36 &   & 409/121 &   \\
14 &   &   & 289/118 & 289/118 &   &   &   \\
14 &   &   &   & 409/121 &   &   &   \\
\hline
15 &   &   & 151/20 & 441/101 & 151/20 & 97/13 & 97/13 \\
15 &   &   & 319/144 & 463/130 & 319/144 & 361/78 & 463/130 \\
15 &   &   & 359/82 &   & 359/82 & 463/130 &   \\
15 &   &   & 463/130 &   & 433/122 &   &   \\
15 &   &   &   &   & 447/98 &   &   \\
15 &   &   &   &   & 463/130 &   &   \\
\hline
16 &   &   & 529/114 & 593/130 & 593/130 & 529/114 & 593/130 \\
16 &   &   & 593/130 &   & 817/239 & 559/122 & 817/239 \\
16 &   &   & 777/208 &   &   & 593/130 &   \\
16 &   &   &   &   &   & 815/237 &   \\
\hline
17 &   &   & 975/274 &   & 703/131 & 1353/380 & 31/2 \\
17 &   &   & 983/260 &   & 1321/288 &   & 127/15 \\
17 &   &   & 1351/362 &   &   &   & 975/274 \\
17 &   &   &   &   &   &   & 983/260 \\

\hline

\end{tabular}
\end{center}
\caption{The sets $L(2, n_y)$ for $3\le n_y\le 15$ given by crossing number and 2-bridge fraction $p/q$.}

\label{p/q for n_y up to 15 part 1}
\end{table}%

\begin{table}[htdp]
\begin{center}
\begin{tabular}{|l|c|c|c|c|c|c|c|}
\hline
&$\L(2,3)$&$\L(2,5)$&$\L(2,7)$&$\L(2,9)$&$\L(2,11)$&$\L(2,13)$&$\L(2,15)$\\
\hline
cr&$p/q$&$p/q$&$p/q$&$p/q$&$p/q$&$p/q$&$p/q$\\
\hline
18 &   &   & 1681/450 & 33/2 & 33/2 & 1489/337 & 241/32 \\
18 &   &   &   & 129/52 & 129/52 &   & 1487/335 \\
18 &   &   &   & 529/214 & 241/32 &   &   \\
18 &   &   &   & 1681/696 & 529/214 &   &   \\
18 &   &   &   & 2321/622 & 1681/696 &   &   \\
\hline
19 &   &   & 2911/780 & 273/32 & 4063/1202 & 273/32 &   \\
19 &   &   &   & 673/78 &   & 673/78 &   \\
19 &   &   &   & 1961/800 &   & 1025/161 &   \\
19 &   &   &   & 2001/898 &   & 1961/800 &   \\
19 &   &   &   & 3329/989 &   & 2001/898 &   \\
19 &   &   &   &   &   & 4015/1106 &   \\
\hline
20 &   &   &   & 3761/1056 & 1279/282 & 3409/744 & 2481/559 \\
20 &   &   &   &   & 3535/996 & 3521/992 & 3631/796 \\
20 &   &   &   &   & 3761/1056 & 3761/1056 & 3761/1056 \\
20 &   &   &   &   & 5681/1661 &   &   \\
\hline
21 &   &   &   & 4297/926 & 4817/1056 & 4817/1056 & 1407/191 \\
21 &   &   &   & 4305/944 & 6143/1407 & 9833/2909 & 4817/1056 \\
21 &   &   &   & 4817/1056 &   &   &   \\
\hline
22 &   &   &   & 7921/2224 & 41/2 & 41/2 & 10737/2354 \\
22 &   &   &   & 7985/2112 & 169/68 & 169/68 &   \\
22 &   &   &   & 10865/2912 & 769/310 & 769/310 &   \\
22 &   &   &   &   & 3073/1272 & 3073/1272 &   \\
22 &   &   &   &   & 9801/4058 & 9801/4058 &   \\
22 &   &   &   &   & 11257/3102 & 10471/3098 &   \\
22 &   &   &   &   &   & 11441/3327 &   \\
\hline
23 &   &   &   & 18817/5042 & 415/36 & 2129/282 & 415/36 \\
23 &   &   &   &   & 1223/146 & 9793/1825 & 1223/146 \\
23 &   &   &   &   & 3487/1424 &   & 3487/1424 \\
23 &   &   &   &   & 11759/4802 &   & 11759/4802 \\
23 &   &   &   &   & 11999/4900 &   & 11999/4900 \\
23 &   &   &   &   & 18663/4996 &   &   \\
\hline
24 &   &   &   & 23409/6272 & 2425/322 & 33433/9892 & 21009/4733 \\
24 &   &   &   &   & 6041/800 &   & 31759/6924 \\
24 &   &   &   &   & 12769/5762 &   &   \\
24 &   &   &   &   & 14369/3282 &   &   \\
24 &   &   &   &   & 32593/8738 &   &   \\
\hline
25 &   &   &   & 40545/10864 & 10359/2284 & 10377/2288 & 14279/2243 \\
25 &   &   &   &   & 16511/3554 & 28743/8098 & 30551/8578 \\
25 &   &   &   &   & 18583/4000 & 30551/8578 &   \\
25 &   &   &   &   & 23543/5086 & 55801/14954 &   \\
25 &   &   &   &   & 30551/8578 &   &   \\
25 &   &   &   &   & 46367/13775 &   &   \\
\hline

\end{tabular}
\end{center}
\caption{The sets $L(2, n_y)$ for $3\le n_y\le 15$ given by crossing number and 2-bridge fraction $p/q$.}

\label{p/q for n_y up to 15 part 2}
\end{table}

\begin{table}[htdp]
\begin{center}
\begin{tabular}{|l|c|c|c|c|c|c|c|}
\hline
&$\L(2,3)$&$\L(2,5)$&$\L(2,7)$&$\L(2,9)$&$\L(2,11)$&$\L(2,13)$&$\L(2,15)$\\
\hline
cr&$p/q$&$p/q$&$p/q$&$p/q$&$p/q$&$p/q$&$p/q$\\

\hline
26 &   &   &   &   & 34905/7522 & 49/2 & 49/2 \\
26 &   &   &   &   & 34969/7666 & 209/84 & 209/84 \\
26 &   &   &   &   & 39129/8578 & 1009/406 & 1009/406 \\
26 &   &   &   &   & 51193/14384 & 4465/1848 & 4465/1848 \\
26 &   &   &   &   &   & 17921/7418 & 17921/7418 \\
26 &   &   &   &   &   & 39129/8578 & 39129/8578 \\
26 &   &   &   &   &   & 57121/23660 & 57121/23660 \\
26 &   &   &   &   &   & 80305/23857 & 97681/28898 \\
\hline
27 &   &   &   &   & 64343/18066 & 601/48 & 3871/514 \\
27 &   &   &   &   & 64351/18068 & 1793/142 &   \\
27 &   &   &   &   & 64863/17156 & 6409/2592 &   \\
27 &   &   &   &   & 87031/23298 & 21401/8738 &   \\
27 &   &   &   &   &   & 69121/28222 &   \\
27 &   &   &   &   &   & 69361/28320 &   \\
27 &   &   &   &   &   & 85561/19597 &   \\
27 &   &   &   &   &   & 91193/25584 &   \\
\hline
28 &   &   &   &   & 108241/28622 & 4385/514 & 162521/48279 \\
28 &   &   &   &   & 108657/28732 & 10817/1268 &   \\
28 &   &   &   &   & 151313/40544 & 26937/3122 &   \\
28 &   &   &   &   &   & 78489/32020 &   \\
28 &   &   &   &   &   & 80089/35940 &   \\
28 &   &   &   &   &   & 151697/44882 &   \\
\hline
29 &   &   &   &   & 188175/50374 & 21441/2840 & 136399/25419 \\
29 &   &   &   &   & 188287/50404 & 50969/11642 &   \\
29 &   &   &   &   & 262087/70226 & 81025/36576 &   \\
29 &   &   &   &   &   & 91193/20830 &   \\
29 &   &   &   &   &   & 115561/26396 &   \\
29 &   &   &   &   &   & 259969/69676 &   \\
\hline
30 &   &   &   &   & 326041/87362 & 84041/18530 & 57/2 \\
30 &   &   &   &   &   & 134689/28992 & 249/100 \\
30 &   &   &   &   &   & 151553/32622 & 1249/502 \\
30 &   &   &   &   &   & 191561/41968 & 5857/2424 \\
30 &   &   &   &   &   & 248169/69680 & 26041/10778 \\
30 &   &   &   &   &   & 453929/121630 & 104441/43260 \\
30 &   &   &   &   &   &   & 248169/69680 \\
30 &   &   &   &   &   &   & 332929/137902 \\
30 &   &   &   &   &   &   & 448689/120112 \\
\hline

\end{tabular}
\end{center}
\caption{The sets $L(2, n_y)$ for $3\le n_y\le 15$ given by crossing number and 2-bridge fraction $p/q$.}

\label{p/q for n_y up to 15 part 3}
\end{table}

\begin{table}[htdp]
\begin{center}
\begin{tabular}{|l|c|c|c|c|c|c|c|}
\hline
&$\L(2,3)$&$\L(2,5)$&$\L(2,7)$&$\L(2,9)$&$\L(2,11)$&$\L(2,13)$&$\L(2,15)$\\
\hline
cr&$p/q$&$p/q$&$p/q$&$p/q$&$p/q$&$p/q$&$p/q$\\

\hline
31 &   &   &   &   & 564719/151316 & 283537/61102 & 807/52 \\
31 &   &   &   &   &   & 284049/62270 & 2599/210 \\
31 &   &   &   &   &   & 284057/62272 & 9215/3728 \\
31 &   &   &   &   &   & 317849/69680 & 39159/15842 \\
31 &   &   &   &   &   & 415929/116866 & 127399/51540 \\
31 &   &   &   &   &   & 645809/191861 & 317849/69680 \\
31 &   &   &   &   &   &   & 400895/166464 \\
31 &   &   &   &   &   &   & 402287/164258 \\
\hline
32 &   &   &   &   &   & 522665/146752 & 1118489/332261 \\
32 &   &   &   &   &   & 522729/146768 &   \\
32 &   &   &   &   &   & 526889/139360 &   \\
32 &   &   &   &   &   & 697745/186784 &   \\
\hline
33 &   &   &   &   &   & 864945/228716 & 38951/5172 \\
33 &   &   &   &   &   & 872153/230622 & 97015/12882 \\
33 &   &   &   &   &   & 882809/233440 & 241791/32020 \\
33 &   &   &   &   &   & 1216977/326170 & 511079/230624 \\
33 &   &   &   &   &   &   & 511119/230642 \\
33 &   &   &   &   &   &   & 575119/131362 \\
33 &   &   &   &   &   &   & 1191711/272951 \\
\hline
34 &   &   &   &   &   & 1509537/404098 & 2090425/559602 \\
34 &   &   &   &   &   & 1515361/407460 &   \\
34 &   &   &   &   &   & 2107561/564720 &   \\
\hline
35 &   &   &   &   &   & 2621905/702714 & 2015903/566018 \\
35 &   &   &   &   &   & 2622017/702520 &   \\
35 &   &   &   &   &   & 3650401/978122 &   \\
\hline
36 &   &   &   &   &   & 4541161/1216800 & 2303201/496338 \\
36 &   &   &   &   &   &   & 2307361/505826 \\
36 &   &   &   &   &   &   & 2307425/505842 \\
36 &   &   &   &   &   &   & 2581921/566018 \\
36 &   &   &   &   &   &   & 6322681/1694162 \\
\hline
37 &   &   &   &   &   & 7865521/2107560 & 4245663/1192082 \\
37 &   &   &   &   &   &   & 4246175/1192226 \\
37 &   &   &   &   &   &   & 4246183/1192212 \\
37 &   &   &   &   &   &   & 4279975/1132036 \\
37 &   &   &   &   &   &   & 8994959/2672279 \\
\hline
38 &   &   &   &   &   &   & 9748249/2609584 \\
\hline
40 &   &   &   &   &   &   & 21077281/5642338 \\
40 &   &   &   &   &   &   & 21089825/5645698 \\
40 &   &   &   &   &   &   & 29354521/7865520 \\
\hline
41 &   &   &   &   &   &   & 50843527/13623482 \\
\hline
42 &   &   &   &   &   &   & 63250209/16947842 \\
\hline
43 &   &   &   &   &   &   & 109552575/29354524 \\
\hline

\end{tabular}
\end{center}
\caption{The sets $L(2, n_y)$ for $3\le n_y\le 15$ given by crossing number and 2-bridge fraction $p/q$.}

\label{p/q for n_y up to 15 part 4}
\end{table}

\begin{table}[htdp]
\begin{center}
\begin{tabular}{|l|c|cccccc|}
\hline
knot&p/q & $n_x$ & $n_y$ & $n_z$ & $\phi_x$ & $\phi_y$ & $\phi_z$\\
\hline
5a1 & $7/2$ & 2 & 3 & 11 & 0 & 0.56099 & 2.58059 \\
\hline
6a3 & $9/2$ & 2 & 3 & 11 & 0 & 0.67319 & 0.89759 \\
\hline
7a3 & $17/5$ & 2 & 5 & 17 & 0 & 0.49979 & 2.64179 \\
7a6 & $15/4$ & 2 & 3 & 11 & 0 & 0.78539 & 2.35619 \\
\hline
8a1 & $31/12$ & 2 & 11 & 41 & 0 & 0.39269 & 2.74889 \\
\hline
9a8 & $31/11$ & 2 & 11 & 41 & 0 & 0.48332 & 1.08747 \\
9a27 & $15/2$ & 2 & 7 & 25 & 0 & 0.49087 & 1.07992 \\
9a33 & $31/7$ & 2 & 7 & 23 & 0 & 0.47123 & 2.67035 \\
\hline
10a23 & $49/20$ & 2 & 5 & 17 & 0 & 0.71399 & 0.85679 \\
10a63 & $55/12$ & 2 & 7 & 25 & 0 & 0.44178 & 2.69980 \\
10a69 & $57/16$ & 2 & 5 & 19 & 0 & 0.58904 & 2.55254 \\
10a75 & $17/2$ & 2 & 5 & 17 & 0 & 0.57119 & 0.99959 \\
\hline
11a91 & $129/49$ & 2 & 41 & 153 & 0 & 0.34816 & 2.79342 \\
11a140 & $65/17$ & 2 & 17 & 63 & 0 & 0.47123 & 1.09955 \\
11a192 & $97/26$ & 2 & 5 & 17 & 0 & 0.64259 & 2.49899 \\
11a210 & $73/16$ & 2 & 5 & 19 & 0 & 0.65449 & 0.91629 \\
11a226 & $71/20$ & 2 & 13 & 49 & 0 & 0.45603 & 1.11475 \\
11a246 & $41/13$ & 2 & 13 & 47 & 0 & 0.47123 & 1.09955 \\
11a333 & $65/14$ & 2 & 5 & 19 & 0 & 0.78539 & 0.78539 \\
11a334 & $49/9$ & 2 & 9 & 29 & 0 & 0.45470 & 2.68688 \\
\hline
12a38 & $71/28$ & 2 & 25 & 89 & 0 & 0.41336 & 1.15742 \\
12a257 & $191/74$ & 2 & 17 & 63 & 0 & 0.37306 & 2.76852 \\
12a715 & $169/50$ & 2 & 7 & 25 & 0 & 0.53996 & 2.60163 \\
12a729 & $167/46$ & 2 & 9 & 31 & 0 & 0.43196 & 2.70962 \\
12a1034 & $121/32$ & 2 & 5 & 19 & 0 & 0.71994 & 2.42164 \\
\hline
13a640 & $55/19$ & 2 & 19 & 65 & 0 & 0.44879 & 1.12199 \\
13a1884 & $289/80$ & 2 & 25 & 93 & 0 & 0.35941 & 2.78217 \\
13a2683 & $287/79$ & 2 & 63 & 235 & 0 & 0.34262 & 2.79896 \\
13a2760 & $239/71$ & 2 & 7 & 23 & 0 & 0.57595 & 2.56563 \\
13a3143 & $23/2$ & 2 & 11 & 37 & 0 & 0.45814 & 1.11264 \\
13a3896 & $111/23$ & 2 & 23 & 85 & 0 & 0.46542 & 1.10537 \\
13a4304 & $209/56$ & 2 & 5 & 17 & 0 & 0.78539 & 2.35619 \\
13a4570 & $79/19$ & 2 & 19 & 69 & 0 & 0.46409 & 1.10669 \\
13a4822 & $71/11$ & 2 & 11 & 35 & 0 & 0.44392 & 2.69767 \\
\hline
\end{tabular}
\end{center}
\caption{All Lissajous knots with frequencies $n_x=2$,  $3\le n_y\le 105$ and with less than 14 crossings. Knot names are as in {\sl Knotscape}.}
\label{2-bridge Lissajous knots 3 to 16 part 1}
\end{table}

\begin{table}[htdp]
\begin{center}
\begin{tabular}{|l|c|cccccc|}
\hline
knot&p/q & $n_x$ & $n_y$ & $n_z$ & $\phi_x$ & $\phi_y$ & $\phi_z$\\
\hline
14a2651 & $89/36$ & 2 & 7 & 23 & 0 & 0.62831 & 0.94247 \\
14a6166 & $289/118$ & 2 & 7 & 23 & 0 & 0.73303 & 0.83775 \\
14a12186 & $407/119$ & 2 & 11 & 37 & 0 & 0.42542 & 2.71616 \\
14a12212 & $409/121$ & 2 & 9 & 31 & 0 & 0.51050 & 2.63108 \\
14a12308 & $103/12$ & 2 & 57 & 215 & 0 & 0.36382 & 2.77776 \\
14a12652 & $127/28$ & 2 & 29 & 103 & 0 & 0.40459 & 1.16619 \\
14a12741 & $25/2$ & 2 & 7 & 23 & 0 & 0.52359 & 1.04719 \\
\hline
15a21965 & $113/29$ & 2 & 29 & 99 & 0 & 0.44178 & 1.12900 \\
15a25723 & $745/288$ & 2 & 19 & 71 & 0 & 0.40142 & 2.74016 \\
15a32142 & $319/144$ & 2 & 7 & 25 & 0 & 0.78539 & 0.78539 \\
15a44612 & $359/82$ & 2 & 7 & 25 & 0 & 0.68722 & 0.88357 \\
15a46260 & $361/78$ & 2 & 13 & 49 & 0 & 0.50670 & 1.06408 \\
15a50643 & $447/98$ & 2 & 11 & 39 & 0 & 0.40840 & 2.73318 \\
15a50772 & $433/122$ & 2 & 11 & 41 & 0 & 0.45311 & 2.68847 \\
15a51438 & $463/130$ & 2 & 7 & 27 & 0 & 0.60059 & 2.54099 \\
15a52567 & $151/20$ & 2 & 7 & 25 & 0 & 0.58904 & 0.98174 \\
15a54893 & $65/21$ & 2 & 21 & 71 & 0 & 0.44392 & 1.12687 \\
15a71359 & $505/109$ & 2 & 85 & 317 & 0 & 0.33994 & 2.80164 \\
15a71603 & $169/29$ & 2 & 29 & 107 & 0 & 0.46199 & 1.10879 \\
15a76044 & $441/101$ & 2 & 9 & 29 & 0 & 0.53737 & 2.60421 \\
15a78853 & $129/25$ & 2 & 25 & 91 & 0 & 0.46040 & 1.11039 \\
15a84772 & $97/13$ & 2 & 13 & 41 & 0 & 0.43633 & 2.70526 \\
\hline
16a7016 & $111/44$ & 2 & 39 & 131 & 0 & 0.40655 & 1.16423 \\
16a57423 & $431/170$ & 2 & 39 & 139 & 0 & 0.40593 & 1.16486 \\
16a135506 & $1103/456$ & 2 & 23 & 85 & 0 & 0.36361 & 2.77798 \\
16a219884 & $961/208$ & 2 & 93 & 347 & 0 & 0.33914 & 2.80244 \\
16a221291 & $777/208$ & 2 & 7 & 25 & 0 & 0.63813 & 2.50345 \\
16a221836 & $783/220$ & 2 & 41 & 155 & 0 & 0.37667 & 2.76492 \\
16a224238 & $791/212$ & 2 & 17 & 65 & 0 & 0.44058 & 2.70100 \\
16a225074 & $593/130$ & 2 & 7 & 27 & 0 & 0.64679 & 0.92399 \\
16a228722 & $559/122$ & 2 & 13 & 47 & 0 & 0.39269 & 2.74889 \\
16a229409 & $577/162$ & 2 & 21 & 79 & 0 & 0.43982 & 1.13097 \\
16a249132 & $817/239$ & 2 & 11 & 37 & 0 & 0.49087 & 2.65071 \\
16a249195 & $815/237$ & 2 & 13 & 43 & 0 & 0.42074 & 2.72084 \\
16a252419 & $385/114$ & 2 & 43 & 153 & 0 & 0.40071 & 1.17008 \\
16a252465 & $399/110$ & 2 & 47 & 177 & 0 & 0.37166 & 2.76993 \\
16a333209 & $529/114$ & 2 & 7 & 27 & 0 & 0.73919 & 0.83159 \\
\hline
\end{tabular}
\end{center}
\caption{All Lissajous knots with frequencies $n_x=2$,  $3\le n_y\le 105$ and with 14--16 crossings. Knot names are as in {\sl Knotscape}.}
\label{2-bridge Lissajous knots 3 to 16 part 2}
\end{table}%

\begin{table}[htdp]
\begin{center}
\begin{tabular}{|l|cccccc|}
\hline
knot& $n_x$ & $n_y$ & $n_z$ & $\phi_x$ & $\phi_y$ & $\phi_z$\\
\hline
5a1 & 3 & 4 & 29 & 0 & 0.49186 & 1.60252\\
\hline
6a3 & 3 & 4 & 29 & 0 & 0.39666 & 1.69772\\
\hline
7a6 & 3 & 7 & 68 & 0 & 0.52359 & 2.61799\\
\hline
8a2 & 3 & 4 & 23 & 0 & 0.29088 & 2.85070\\
8n2 & 3 & 4 & 37 & 0 & 0.49805 & 2.64353\\
\hline
9a25 & 3 & 5 & 28 & 0 & 0.26973 & 1.82466\\
\hline
10a20 & 3 & 4 & 23 & 0 & 0.32967 & 0.71752\\
10a73 & 3 & 7 & 40 & 0 & 0.23394 & 2.90764\\
10a89 & 3 & 4 & 29 & 0 & 0.36493 & 0.68226\\
10a103 & 3 & 5 & 29 & 0 & 0.41579 & 2.72579\\
10n28 & 3 & 7 & 40 & 0 & 0.21166 & 1.88272\\
\hline
11n50 & 4 & 5 & 57 & 0 & 0.36736 & 1.98883\\
11n151 & 4 & 7 & 69 & 0 & 0.24802 & 1.32277\\
\hline
12a426 & 3 & 10 & 83 & 0 & 0.45603 & 2.68555\\
12a427 & 3 & 5 & 29 & 0 & 0.32339 & 2.81819\\
12a448 & 3 & 5 & 28 & 0 & 0.30146 & 0.74573\\
12a868 & 3 & 5 & 26 & 0 & 0.28713 & 1.80726\\
12a1164 & 4 & 5 & 61 & 0 & 0.34509 & 2.01109\\
12n133 & 3 & 4 & 23 & 0 & 0.36845 & 1.72593\\
12n293 & 3 & 10 & 59 & 0 & 0.28077 & 1.81362\\
12n322 & 4 & 5 & 37 & 0 & 0.16829 & 2.18789\\
12n483 & 5 & 6 & 113 & 0 & 0.35111 & 2.16215\\
\hline
13a2233 & 3 & 8 & 67 & 0 & 0.48171 & 2.65988\\
13a4774 & 3 & 4 & 23 & 0 & 0.40724 & 2.73434\\
13n1405 & 4 & 5 & 51 & 0 & 0.35062 & 2.00557\\
13n1734 & 3 & 10 & 103 & 0 & 0.56066 & 1.53372\\
13n3594 & 3 & 4 & 23 & 0 & 0.44602 & 0.60116\\
\hline
\end{tabular}
\end{center}
\caption{Small-crossing Lissajous knots with all frequencies greater than 2. Only three knots, $5a1$, $6a3$ and $7a6$ are 2-bridge and appear in Table~\ref{2-bridge Lissajous knots 3 to 16 part 1}. Knot names are as in {\sl Knotscape}.}
\label{Lissajous knots part 1}
\end{table}%

\begin{table}[htdp]
\begin{center}
\begin{tabular}{|l|cccccc|}
\hline
knot& $n_x$ & $n_y$ & $n_z$ & $\phi_x$ & $\phi_y$ & $\phi_z$\\
\hline
14a1491 & 4 & 5 & 39 & 0 & 0.19634 & 0.58904\\
14a6398 & 3 & 7 & 59 & 0 & 0.49979 & 2.64179\\
14a6912 & 3 & 8 & 67 & 0 & 0.42586 & 0.62133\\
14a8662 & 3 & 7 & 53 & 0 & 0.35779 & 0.68940\\
14a13089 & 4 & 7 & 53 & 0 & 0.15707 & 1.41371\\
14a15296 & 3 & 4 & 29 & 0 & 0.46013 & 0.58706\\
14a16309 & 3 & 7 & 53 & 0 & 0.47996 & 1.61442\\
14a16437 & 3 & 5 & 29 & 0 & 0.35419 & 1.74019\\
14a18187 & 3 & 5 & 37 & 0 & 0.38646 & 0.66073\\
14n6560 & 3 & 8 & 47 & 0 & 0.25703 & 2.88455\\
14n9732 & 3 & 7 & 53 & 0 & 0.42760 & 1.66678\\
14n13886 & 3 & 4 & 29 & 0 & 0.42839 & 2.71319\\
14n14189 & 4 & 5 & 63 & 0 & 0.39269 & 2.74889\\
14n15552 & 3 & 5 & 28 & 0 & 0.33319 & 2.80839\\
14n18513 & 3 & 8 & 83 & 0 & 0.58113 & 1.51325\\
14n22071 & 3 & 7 & 64 & 0 & 0.50884 & 2.63274\\
14n22073 & 3 & 5 & 29 & 0 & 0.50819 & 2.63339\\
14n23738 & 3 & 7 & 50 & 0 & 0.39499 & 1.69939\\
14n24494 & 3 & 8 & 61 & 0 & 0.35665 & 1.73774\\
14n25903 & 3 & 7 & 53 & 0 & 0.37524 & 1.71914\\
\hline
15a80928 & 4 & 5 & 37 & 0 & 0.20569 & 0.57969\\
15n77228 & 4 & 5 & 69 & 0 & 0.42453 & 1.14625\\
15n92508 & 3 & 8 & 61 & 0 & 0.43253 & 2.70905\\
15n103019 & 5 & 6 & 103 & 0 & 0.30263 & 2.83896\\
15n116110 & 4 & 7 & 53 & 0 & 0.14398 & 2.21220\\
\hline
16a128851 & 3 & 7 & 40 & 0 & 0.30079 & 2.84080\\
16a151023 & 3 & 7 & 38 & 0 & 0.26761 & 0.77958\\
16a168328 & 3 & 7 & 41 & 0 & 0.29452 & 2.84706\\
16a202258 & 3 & 7 & 40 & 0 & 0.27850 & 1.81588\\
16a295212 & 3 & 5 & 26 & 0 & 0.38847 & 1.70591\\
16a312423 & 3 & 5 & 29 & 0 & 0.38499 & 0.66219\\
16a340770 & 3 & 7 & 41 & 0 & 0.38179 & 0.66540\\
16n42863 & 4 & 5 & 39 & 0 & 0.21419 & 1.35659\\
16n228473 & 3 & 7 & 38 & 0 & 0.29088 & 1.80350\\
16n390014 & 3 & 5 & 26 & 0 & 0.35469 & 2.78689\\
16n507235 & 3 & 5 & 28 & 0 & 0.39666 & 0.65053\\
16n562396 & 3 & 5 & 28 & 0 & 0.36493 & 1.72946\\
16n768985 & 3 & 4 & 23 & 0 & 0.48481 & 1.60958\\
16n982564 & 4 & 5 & 73 & 0 & 0.44304 & 1.12775\\
16n988939 & 3 & 7 & 41 & 0 & 0.51269 & 0.53450\\
16n1008347 & 3 & 7 & 50 & 0 & 0.32150 & 0.72568\\
\hline

\end{tabular}
\end{center}
\caption{Small-crossing Lissajous knots with all frequencies greater than 2.  Knot names are as in {\sl Knotscape}.}
\label{Lissajous knots part 2}
\end{table}%

\begin{table}[htdp]
\begin{center}
\begin{tabular}{|l|c|cccccccc|}
\hline
knot& $p/q$&$n_x$ & $n_y$ & $n_{z,1}$&$n_{z,2}$ & $\phi_x$ & $\phi_y$ & $\phi_{z,1}$&$\phi_{z,2}$\\
\hline
3a1 & 3$/$1 & 2 & 3 & 1 & 2 & 0 & $\pi/4$ & 0.39269 & 1.66017 \\
\hline
4a1 & 5$/$2 & 2 & 3 & 1 & 3 & 0 & $\pi/4$ & 1.62773 & 5.79254 \\
\hline
{\bf 5a1} & {\bf 7$/$2} & 2 & 3 & 1 & 7 & 0 & $\pi/4$ & 0.10580 & 2.49320 \\
5a2 & 5$/$1 & 2 & 5 & 2 & 3 & 0 & $\pi/4$ & 0.96046 & 4.09767 \\
\hline
6a1 & 13$/$5 & 2 & 5 & 1 & 5 & 0 & $\pi/4$ & 0.03573 & 2.53353 \\
6a2 & 11$/$3 & 2 & 7 & 1 & 7 & 0 & $\pi/4$ & 1.90655 & 5.01637 \\
{\bf 6a3} & {\bf 9$/$2} & 2 & 3 & 1 & 5 & 0 & $\pi/4$ & 0.18165 & 1.75945 \\
\hline
7a1 & 21$/$8 & 2 & 5 & 1 & 5 & 0 & $\pi/4$ & 1.60021 & 5.52412 \\
7a2 & 19$/$7 & 2 & 7 & 3 & 7 & 0 & $\pi/4$ & 1.66835 & 6.11271 \\
{\bf 7a3} & {\bf 17$/$5} & 2 & 7 & 5 & 9 & 0 & $\pi/4$ & 0.08774 & 5.55745 \\
7a4 & 11$/$2 & 2 & 7 & 3 & 7 & 0 & $\pi/4$ & 1.60853 & 6.27384 \\
7a5 & 13$/$3 & 2 & 9 & 4 & 7 & 0 & $\pi/4$ & 1.57817 & 4.41032 \\
{\bf 7a6} & {\bf 15$/$4} & 2 & 3 & 1 & 7 & 0 & $\pi/4$ & 1.87792 & 4.64352 \\
7a7 & 7$/$1 & 2 & 7 & 2 & 5 & 0 & $\pi/4$ & 0.93991 & 0.93104 \\
\hline
{\bf 8a1} & {\bf 31$/$12} & 2 & 11 & 1 & 5 & 0 & $\pi/4$ & 0.37204 & 1.78795 \\
{\it 8a4} & {\it 25$/$9} & 2 & 7 & 1 & 5 & 0 & $\pi/4$ & 2.04720 & 5.29197 \\
8a5 & 29$/$12 & 2 & 5 & 3 & 5 & 0 & $\pi/4$ & 1.59453 & 2.05821 \\
8a6 & 23$/$5 & 2 & 9 & 1 & 9 & 0 & $\pi/4$ & 0.35397 & 2.65710 \\
8a7 & 29$/$8 & 2 & 9 & 1 & 5 & 0 & $\pi/4$ & 1.47451 & 2.10447 \\
8a8 & 17$/$3 & 2 & 11 & 3 & 10 & 0 & $\pi/4$ & 0.39241 & 5.09182 \\
8a9 & 27$/$8 & 2 & 9 & 1 & 5 & 0 & $\pi/4$ & 2.03830 & 2.05668 \\
{\it 8a10} & {\it 23$/$7} & 2 & 7 & 1 & 9 & 0 & $\pi/4$ & 0.48400 & 5.18915 \\
8a11 & 13$/$2 & 2 & 5 & 3 & 5 & 0 & $\pi/4$ & 1.58524 & 0.24531 \\
8a16 & 25$/$7 & 2 & 7 & 3 & 7 & 0 & $\pi/4$ & 0.04412 & 2.25248 \\
8a17 & 19$/$4 & 2 & 9 & 3 & 7 & 0 & $\pi/4$ & 1.92077 & 6.06457 \\
{\it 8a18} & {\it 17$/$4} & 2 & 7 & 1 & 5 & 0 & $\pi/4$ & 1.42912 & 1.98797 \\
\hline
{\it 9a3} & {\it 41$/$16} & 2 & 7 & 3 & 5 & 0 & $\pi/4$ & 1.88608 & 4.98854 \\
{\bf 9a8} & {\bf 31$/$11} & 2 & 5 & 7 & 9 & 0 & $\pi/4$ & 0.09293 & 5.50516 \\
{\it 9a10} & {\it 39$/$16} & 2 & 5 & 1 & 7 & 0 & $\pi/4$ & 2.08636 & 5.17367 \\
9a12 & 49$/$18 & 2 & 17 & 1 & 9 & 0 & $\pi/4$ & 0.60289 & 4.56332 \\
9a13 & 55$/$21 & 2 & 7 & 1 & 7 & 0 & $\pi/4$ & 0.00613 & 5.43134 \\
9a14 & 39$/$14 & 2 & 9 & 1 & 9 & 0 & $\pi/4$ & 2.13083 & 5.89035 \\
9a15 & 47$/$13 & 2 & 9 & 1 & 9 & 0 & $\pi/4$ & 1.74912 & 1.92322 \\
9a16 & 45$/$19 & 2 & 9 & 1 & 5 & 0 & $\pi/4$ & 1.93719 & 4.94328 \\
9a17 & 37$/$8 & 2 & 5 & 2 & 9 & 0 & $\pi/4$ & 0.19367 & 2.96295 \\
9a19 & 41$/$11 & 2 & 15 & 3 & 10 & 0 & $\pi/4$ & 2.62089 & 0.60844 \\
9a20 & 33$/$7 & 2 & 17 & 5 & 7 & 0 & $\pi/4$ & 0.10752 & 5.13086 \\
9a21 & 43$/$12 & 2 & 5 & 3 & 9 & 0 & $\pi/4$ & 2.14045 & 5.61205 \\
9a22 & 35$/$8 & 2 & 5 & 7 & 9 & 0 & $\pi/4$ & 1.67486 & 1.65979 \\
9a23 & 27$/$5 & 2 & 13 & 4 & 5 & 0 & $\pi/4$ & 1.36285 & 5.38881 \\
{\it 9a24} & {\it 41$/$12} & 2 & 11 & 1 & 7 & 0 & $\pi/4$ & 0.44828 & 2.24339 \\
9a26 & 29$/$9 & 2 & 11 & 8 & 9 & 0 & $\pi/4$ & 1.42268 & 3.50098 \\
{\bf 9a27} & {\bf 15$/$2} & 2 & 5 & 1 & 7 & 0 & $\pi/4$ & 0.00135 & 6.21184 \\
{\bf 9a33} & {\bf 31$/$7} & 2 & 11 & 4 & 9 & 0 & $\pi/4$ & 0.93417 & 3.86038 \\
9a34 & 37$/$10 & 2 & 9 & 3 & 5 & 0 & $\pi/4$ & 1.98367 & 5.56618 \\
9a35 & 21$/$4 & 2 & 11 & 9 & 10 & 0 & $\pi/4$ & 0.35932 & 5.18305 \\
{\it 9a36} & {\it 23$/$4} & 2 & 9 & 3 & 5 & 0 & $\pi/4$ & 1.65102 & 3.04593 \\
9a38 & 19$/$3 & 2 & 13 & 1 & 7 & 0 & $\pi/4$ & 1.93386 & 2.02910 \\
{\it 9a39} & {\it 33$/$10} & 2 & 11 & 3 & 7 & 0 & $\pi/4$ & 2.16159 & 2.03213 \\
9a41 & 9$/$1 & 2 & 9 & 2 & 7 & 0 & $\pi/4$ & 0.86114 & 0.74621 \\

\hline
\end{tabular}
\end{center}
\caption{Fourier-$(1,1,2)$ descriptions of all 2-bridge knots to 9 crossings. All amplitudes are 1. Boldface entries are  known to be Lissajous while italic entries {\it might} be Lissajous. All others cannot be Lissajous. Knot names are as in {\sl Knotscape}.}
\label{2-bridge 3 to 10 as fourier(1,1,2), part 1}
\end{table}%

\begin{table}[htdp]
\begin{center}
\begin{tabular}{|l|c|cccccccc|}
\hline
knot& $p/q$&$n_x$ & $n_y$ & $n_{z,1}$&$n_{z,2}$ & $\phi_x$ & $\phi_y$ & $\phi_{z,1}$&$\phi_{z,2}$\\
\hline
10a5 & 51$/$20 & 2 & 11 & 3 & 7 & 0 & $\pi/4$ & 0.18595 & 2.82470 \\
10a19 & 37$/$13 & 2 & 13 & 3 & 7 & 0 & $\pi/4$ & 1.85642 & 1.83386 \\
{\bf 10a23} & {\bf 49$/$20} & 2 & 5 & 1 & 7 & 0 & $\pi/4$ & 0.36657 & 2.47281 \\
10a25 & 89$/$34 & 2 & 7 & 1 & 7 & 0 & $\pi/4$ & 1.62228 & 5.59775 \\
10a26 & 61$/$22 & 2 & 17 & 3 & 8 & 0 & $\pi/4$ & 0.19982 & 4.43387 \\
10a29 & 59$/$25 & 2 & 9 & 1 & 5 & 0 & $\pi/4$ & 1.82942 & 5.29828 \\
10a30 & 75$/$29 & 2 & 15 & 5 & 7 & 0 & $\pi/4$ & 1.66311 & 4.79250 \\
10a31 & 81$/$31 & 2 & 17 & 3 & 7 & 0 & $\pi/4$ & 1.36514 & 5.07457 \\
10a32 & 79$/$29 & 2 & 17 & 3 & 7 & 0 & $\pi/4$ & 2.51922 & 6.10032 \\
{\it 10a33} & {\it 57$/$13} & 2 & 17 & 1 & 9 & 0 & $\pi/4$ & 1.22778 & 5.48698 \\
10a34 & 67$/$18 & 2 & 13 & 1 & 5 & 0 & $\pi/4$ & 1.94658 & 5.20996 \\
10a35 & 71$/$26 & 2 & 17 & 5 & 7 & 0 & $\pi/4$ & 2.17223 & 4.46442 \\
{\it 10a43} & {\it 47$/$11} & 2 & 17 & 5 & 8 & 0 & $\pi/4$ & 2.17425 & 5.50211 \\
10a44 & 53$/$14 & 2 & 9 & 3 & 5 & 0 & $\pi/4$ & 2.12967 & 5.39832 \\
10a49 & 53$/$23 & 2 & 13 & 1 & 7 & 0 & $\pi/4$ & 1.00888 & 3.57446 \\
10a52 & 73$/$27 & 2 & 17 & 3 & 7 & 0 & $\pi/4$ & 2.90021 & 2.80121 \\
10a53 & 63$/$17 & 2 & 17 & 1 & 9 & 0 & $\pi/4$ & 0.46162 & 2.91082 \\
10a54 & 53$/$12 & 2 & 5 & 5 & 7 & 0 & $\pi/4$ & 0.00438 & 1.96518 \\
10a55 & 69$/$19 & 2 & 13 & 3 & 7 & 0 & $\pi/4$ & 2.09996 & 5.23165 \\
10a56 & 33$/$5 & 2 & 19 & 5 & 6 & 0 & $\pi/4$ & 0.37923 & 4.72239 \\
10a57 & 59$/$18 & 2 & 21 & 3 & 8 & 0 & $\pi/4$ & 0.38316 & 4.48117 \\
{\it 10a58} & {\it 71$/$21} & 2 & 13 & 6 & 9 & 0 & $\pi/4$ & 1.04484 & 0.64343 \\
10a59 & 23$/$3 & 2 & 23 & 2 & 7 & 0 & $\pi/4$ & 0.57078 & 0.10280 \\
10a60 & 45$/$14 & 2 & 17 & 1 & 9 & 0 & $\pi/4$ & 1.72208 & 1.81144 \\
{\it 10a61} & {\it 65$/$19} & 2 & 17 & 2 & 5 & 0 & $\pi/4$ & 1.05935 & 4.53128 \\
{\bf 10a63} & {\bf 55$/$12} & 2 & 11 & 1 & 7 & 0 & $\pi/4$ & 2.07884 & 1.80942 \\
10a64 & 45$/$8 & 2 & 17 & 5 & 9 & 0 & $\pi/4$ & 1.41312 & 1.49892 \\
10a65 & 43$/$8 & 2 & 19 & 1 & 7 & 0 & $\pi/4$ & 1.64470 & 4.03256 \\
10a68 & 43$/$9 & 2 & 13 & 9 & 10 & 0 & $\pi/4$ & 0.30693 & 1.31007 \\
{\bf 10a69} & {\bf 57$/$16} & 2 & 9 & 3 & 5 & 0 & $\pi/4$ & 0.02223 & 2.46173 \\
10a70 & 37$/$7 & 2 & 19 & 5 & 7 & 0 & $\pi/4$ & 0.36129 & 4.69287 \\
{\it 10a71} & {\it 55$/$16} & 2 & 13 & 3 & 5 & 0 & $\pi/4$ & 1.15846 & 5.71558 \\
10a74 & 35$/$11 & 2 & 13 & 5 & 10 & 0 & $\pi/4$ & 0.31052 & 4.18157 \\
{\bf 10a75} & {\bf 17$/$2} & 2 & 5 & 2 & 7 & 0 & $\pi/4$ & 0.00295 & 3.10669 \\
10a107 & 41$/$9 & 2 & 9 & 5 & 9 & 0 & $\pi/4$ & 0.11322 & 5.44496 \\
10a108 & 51$/$11 & 2 & 17 & 2 & 7 & 0 & $\pi/4$ & 1.01343 & 4.49878 \\
{\it 10a109} & {\it 65$/$18} & 2 & 11 & 2 & 5 & 0 & $\pi/4$ & 0.22667 & 5.30925 \\
10a110 & 39$/$7 & 2 & 19 & 4 & 5 & 0 & $\pi/4$ & 1.10277 & 4.42691 \\
10a111 & 61$/$17 & 2 & 19 & 4 & 5 & 0 & $\pi/4$ & 1.35067 & 1.23088 \\
{\it 10a112} & {\it 49$/$13} & 2 & 17 & 5 & 7 & 0 & $\pi/4$ & 0.06896 & 1.24086 \\
10a113 & 27$/$4 & 2 & 17 & 4 & 7 & 0 & $\pi/4$ & 1.05187 & 2.20380 \\
10a114 & 29$/$5 & 2 & 19 & 5 & 7 & 0 & $\pi/4$ & 0.18690 & 4.45167 \\
{\it 10a115} & {\it 47$/$10} & 2 & 11 & 7 & 10 & 0 & $\pi/4$ & 0.50632 & 4.53216 \\
10a116 & 43$/$10 & 2 & 15 & 7 & 9 & 0 & $\pi/4$ & 0.09320 & 1.94013 \\
{\it 10a117} & {\it 25$/$4} & 2 & 7 & 1 & 5 & 0 & $\pi/4$ & 2.15907 & 2.20357 \\

\hline
\end{tabular}
\end{center}
\caption{Fourier-$(1,1,2)$ descriptions of all 2-bridge knots with 10 crossings. All amplitudes are 1. Boldface entries are  known to be Lissajous while italic entries {\it might} be Lissajous. All others cannot be Lissajous. Knot names are as in {\sl Knotscape}.}
\label{2-bridge 3 to 10 as fourier(1,1,2), part 2}
\end{table}%

\begin{table}[htdp]
\begin{center}
\begin{tabular}{|l|cccccccc|}
\hline
knot& $n_x$ & $n_y$ & $n_{z,1}$&$n_{z,2}$ & $\phi_x$ & $\phi_y$ & $\phi_{z,1}$&$\phi_{z,2}$\\
\hline
8a2 & 3 & 4 & 1 & 7 & 0 & $\pi/6$ & 1.63362 & 2.03575 \\
8a3 & 3 & 5 & 6 & 13 & 0 & $\pi/6$ & 0.56548 & 2.03575 \\
8a12 & 3 & 4 & 3 & 5 & 0 & $\pi/6$ & 1.04300 & 0.80424 \\
8a13 & 3 & 4 & 1 & 9 & 0 & $\pi/6$ & 0.26389 & 1.58336 \\
8a14 & 3 & 4 & 7 & 14 & 0 & $\pi/6$ & 1.28176 & 1.78442 \\
8a15 & 3 & 7 & 1 & 10 & 0 & $\pi/6$ & 1.64619 & 2.31221 \\
8n1 & 3 & 4 & 1 & 14 & 0 & $\pi/6$ & 1.94778 & 2.76460 \\
8n2 & 3 & 4 & 1 & 5 & 0 & $\pi/6$ & 0.05026 & 2.23681 \\
8n3 & 3 & 4 & 1 & 3 & 0 & $\pi/6$ & 0.26389 & 1.58336 \\
\hline
9a1 & 3 & 5 & 7 & 10 & 0 & $\pi/6$ & 2.29964 & 0.03769 \\
9a2 & 3 & 7 & 1 & 6 & 0 & $\pi/6$ & 0.35185 & 1.05557 \\
9a4 & 3 & 4 & 1 & 14 & 0 & $\pi/6$ & 1.33203 & 2.27451 \\
9a5 & 3 & 8 & 1 & 8 & 0 & $\pi/6$ & 0.65345 & 1.70902 \\
9a6 & 3 & 7 & 1 & 10 & 0 & $\pi/6$ & 1.88495 & 2.43787 \\
9a7 & 3 & 4 & 1 & 14 & 0 & $\pi/6$ & 2.03575 & 2.37504 \\
9a9 & 3 & 7 & 4 & 15 & 0 & $\pi/6$ & 1.15610 & 0.76654 \\
9a11 & 3 & 5 & 9 & 14 & 0 & $\pi/6$ & 1.04300 & 1.33203 \\
9a18 & 3 & 4 & 2 & 7 & 0 & $\pi/6$ & 2.37504 & 2.03575 \\
9a25 & 3 & 7 & 2 & 14 & 0 & $\pi/6$ & 0.98017 & 0.37699 \\
9a28 & 3 & 7 & 4 & 5 & 0 & $\pi/6$ & 1.28176 & 0.45238 \\
9a29 & 4 & 7 & 2 & 13 & 0 & $\pi/8$ & 0.26389 & 2.07345 \\
9a30 & 3 & 7 & 8 & 9 & 0 & $\pi/6$ & 0.23876 & 0.77911 \\
9a31 & 3 & 4 & 10 & 11 & 0 & $\pi/6$ & 1.38230 & 1.87238 \\
9a32 & 3 & 7 & 4 & 13 & 0 & $\pi/6$ & 0.46495 & 1.20637 \\
9a37 & 3 & 4 & 2 & 11 & 0 & $\pi/6$ & 1.04300 & 2.62637 \\
9a40 & 3 & 7 & 4 & 13 & 0 & $\pi/6$ & 1.06814 & 1.06814 \\
9n1 & 3 & 4 & 1 & 14 & 0 & $\pi/6$ & 0.05026 & 0.27646 \\
9n2 & 3 & 5 & 4 & 7 & 0 & $\pi/6$ & 0.15079 & 1.99805 \\
9n3 & 3 & 8 & 1 & 6 & 0 & $\pi/6$ & 0.76654 & 0.95504 \\
9n4 & 3 & 4 & 2 & 11 & 0 & $\pi/6$ & 0.35185 & 2.51327 \\
9n5 & 3 & 4 & 1 & 4 & 0 & $\pi/6$ & 1.33203 & 2.09858 \\
9n6 & 3 & 4 & 2 & 13 & 0 & $\pi/6$ & 0.08796 & 2.48814 \\
9n7 & 3 & 7 & 2 & 9 & 0 & $\pi/6$ & 2.62637 & 1.05557 \\
9n8 & 3 & 7 & 4 & 13 & 0 & $\pi/6$ & 0.05026 & 0.18849 \\
\hline
\end{tabular}
\end{center}
\caption{Fourier-$(1,1,2)$ descriptions of  non 2-bridge knots up to 9 crossings. All amplitudes are 1.  Knot names are as in {\sl Knotscape}.}
\label{non 2-bridge 3 to 10 as fourier(1,1,2), part 1}
\end{table}%

\begin{table}[htdp]
\begin{center}
\begin{tabular}{|l|cccccccc|}
\hline
knot& $n_x$ & $n_y$ & $n_{z,1}$&$n_{z,2}$ & $\phi_x$ & $\phi_y$ & $\phi_{z,1}$&$\phi_{z,2}$\\
\hline
10a1 & 3 & 4 & 5 & 9 & 0 & $\pi/6$ & 1.04300 & 1.58336 \\
10a2 & 3 & 8 & 4 & 13 & 0 & $\pi/6$ & 0.70371 & 0.08796 \\
10a3 & 3 & 8 & 7 & 14 & 0 & $\pi/6$ & 2.70176 & 0.23876 \\
10a4 & 3 & 5 & 4 & 13 & 0 & $\pi/6$ & 2.34991 & 2.34991 \\
10a7 & 3 & 5 & 8 & 11 & 0 & $\pi/6$ & 0.01256 & 1.60849 \\
10a9 & 3 & 8 & 5 & 6 & 0 & $\pi/6$ & 0.27646 & 1.58336 \\
10a10 & 3 & 10 & 2 & 11 & 0 & $\pi/6$ & 2.48814 & 2.94053 \\
10a11 & 3 & 8 & 3 & 14 & 0 & $\pi/6$ & 0.05026 & 1.01787 \\
10a12 & 3 & 8 & 2 & 5 & 0 & $\pi/6$ & 1.40743 & 1.67132 \\
10a14 & 3 & 8 & 2 & 7 & 0 & $\pi/6$ & 0.33929 & 2.81486 \\
10a15 & 4 & 7 & 1 & 12 & 0 & $\pi/8$ & 1.01787 & 0.67858 \\
10a16 & 4 & 7 & 1 & 12 & 0 & $\pi/8$ & 0.20106 & 0.66601 \\
10a17 & 3 & 8 & 5 & 15 & 0 & $\pi/6$ & 0.22619 & 1.06814 \\
10a18 & 3 & 8 & 4 & 9 & 0 & $\pi/6$ & 0.17592 & 1.58336 \\
10a20 & 3 & 4 & 1 & 7 & 0 & $\pi/6$ & 0.22619 & 1.47026 \\
10a21 & 3 & 8 & 1 & 10 & 0 & $\pi/6$ & 0.82938 & 2.71433 \\
10a22 & 3 & 8 & 1 & 10 & 0 & $\pi/6$ & 2.46300 & 2.75203 \\
10a28 & 4 & 7 & 1 & 14 & 0 & $\pi/8$ & 0.43982 & 0.11309 \\
10a36 & 4 & 5 & 1 & 14 & 0 & $\pi/8$ & 0.82938 & 1.60849 \\
10a37 & 3 & 5 & 7 & 14 & 0 & $\pi/6$ & 0.07539 & 2.90283 \\
10a38 & 3 & 7 & 2 & 11 & 0 & $\pi/6$ & 0.66601 & 0.01256 \\
10a42 & 3 & 10 & 5 & 8 & 0 & $\pi/6$ & 2.61380 & 0.42725 \\
10a47 & 3 & 4 & 5 & 14 & 0 & $\pi/6$ & 2.19911 & 1.96035 \\
10a48 & 3 & 8 & 5 & 10 & 0 & $\pi/6$ & 2.02318 & 1.39486 \\
10a50 & 3 & 7 & 1 & 8 & 0 & $\pi/6$ & 1.60849 & 1.09327 \\
10a51 & 3 & 7 & 2 & 4 & 0 & $\pi/6$ & 2.56353 & 5.45380 \\
10a62 & 5 & 6 & 4 & 11 & 0 & $\pi/10$ & 0.27646 & 0.59061 \\
10a66 & 3 & 5 & 2 & 11 & 0 & $\pi/6$ & 1.52053 & 1.85982 \\
10a67 & 3 & 8 & 2 & 13 & 0 & $\pi/6$ & 1.40743 & 3.12902 \\
\hline
\end{tabular}
\end{center}
\caption{Fourier-$(1,1,2)$ descriptions of  alternating non 2-bridge knots with 10 crossings. All amplitudes are 1.  Knot names are as in {\sl Knotscape}.}
\label{non 2-bridge 3 to 10 as fourier(1,1,2), part 2}
\end{table}%

\begin{table}[htdp]
\begin{center}
\begin{tabular}{|l|cccccccc|}
\hline
knot& $n_x$ & $n_y$ & $n_{z,1}$&$n_{z,2}$ & $\phi_x$ & $\phi_y$ & $\phi_{z,1}$&$\phi_{z,2}$\\
\hline
10a72 & 4 & 5 & 7 & 14 & 0 & $\pi/8$ & 2.47557 & 1.73415 \\
10a73 & 3 & 5 & 2 & 10 & 0 & $\pi/6$ & 1.98548 & 0.16336 \\
10a76 & 3 & 7 & 6 & 11 & 0 & $\pi/6$ & 1.58336 & 1.05557 \\
10a77 & 5 & 6 & 3 & 14 & 0 & $\pi/10$ & 0.30159 & 1.75929 \\
10a78 & 5 & 6 & 1 & 12 & 0 & $\pi/10$ & 1.45769 & 1.26920 \\
10a80 & 3 & 10 & 4 & 15 & 0 & $\pi/6$ & 1.97292 & 0.31415 \\
10a82 & 5 & 6 & 4 & 11 & 0 & $\pi/10$ & 0.37699 & 2.86513 \\
10a84 & 4 & 7 & 6 & 13 & 0 & $\pi/8$ & 0.31415 & 1.73415 \\
10a85 & 3 & 4 & 7 & 14 & 0 & $\pi/6$ & 1.08070 & 2.04831 \\
10a87 & 4 & 7 & 2 & 3 & 0 & $\pi/8$ & 1.43256 & 2.14884 \\
10a89 & 3 & 4 & 3 & 13 & 0 & $\pi/6$ & 0.18849 & 1.48283 \\
10a90 & 5 & 6 & 6 & 7 & 0 & $\pi/10$ & 0.02513 & 1.53309 \\
10a91 & 4 & 7 & 1 & 14 & 0 & $\pi/8$ & 0.70371 & 0.05026 \\
10a92 & 3 & 4 & 7 & 14 & 0 & $\pi/6$ & 1.09327 & 2.46300 \\
10a93 & 3 & 8 & 2 & 13 & 0 & $\pi/6$ & 1.57079 & 0.69115 \\
10a94 & 3 & 4 & 7 & 10 & 0 & $\pi/6$ & 1.28176 & 1.36973 \\
10a95 & 4 & 7 & 4 & 13 & 0 & $\pi/8$ & 1.58336 & 0.06283 \\
10a96 & 3 & 4 & 7 & 13 & 0 & $\pi/6$ & 1.63362 & 5.17734 \\
10a97 & 3 & 8 & 2 & 7 & 0 & $\pi/6$ & 1.35716 & 2.36247 \\
10a99 & 3 & 8 & 2 & 7 & 0 & $\pi/6$ & 1.36973 & 2.63893 \\
10a100 & 3 & 10 & 4 & 7 & 0 & $\pi/6$ & 1.73415 & 1.41999 \\
10a101 & 5 & 6 & 3 & 8 & 0 & $\pi/10$ & 0.06283 & 2.57610 \\
10a102 & 3 & 4 & 1 & 13 & 0 & $\pi/6$ & 2.02318 & 1.24407 \\
10a103 & 3 & 5 & 9 & 13 & 0 & $\pi/6$ & 0.75398 & 2.57610 \\
10a105 & 3 & 4 & 10 & 13 & 0 & $\pi/6$ & 1.77185 & 2.86513 \\
10a119 & 4 & 7 & 10 & 15 & 0 & $\pi/8$ & 0.45238 & 1.25663 \\
10a121 & 3 & 5 & 3 & 7 & 0 & $\pi/6$ & 1.04300 & 0.85451 \\
10a122 & 4 & 7 & 7 & 9 & 0 & $\pi/8$ & 1.28176 & 3.99610 \\
10a123 & 3 & 7 & 2 & 6 & 0 & $\pi/6$ & 0.57805 & 4.31026 \\
\hline
\end{tabular}
\end{center}
\caption{Fourier-$(1,1,2)$ descriptions of  alternating non 2-bridge knots with 10 crossings. All amplitudes are 1.  Knot names are as in {\sl Knotscape}.}
\label{non 2-bridge 3 to 10 as fourier(1,1,2), part 3}
\end{table}%

\begin{table}[htdp]
\begin{center}
\begin{tabular}{|l|cccccccc|}
\hline
knot& $n_x$ & $n_y$ & $n_{z,1}$&$n_{z,2}$ & $\phi_x$ & $\phi_y$ & $\phi_{z,1}$&$\phi_{z,2}$\\
\hline

10n1 & 3 & 4 & 3 & 7 & 0 & $\pi/6$ & 1.04300 & 1.47026 \\
10n2 & 3 & 4 & 5 & 14 & 0 & $\pi/6$ & 0.99274 & 2.87769 \\
10n3 & 3 & 4 & 1 & 5 & 0 & $\pi/6$ & 1.94778 & 5.88106 \\
10n4 & 3 & 4 & 6 & 11 & 0 & $\pi/6$ & 0.35185 & 0.18849 \\
10n5 & 3 & 4 & 1 & 10 & 0 & $\pi/6$ & 1.11840 & 2.37504 \\
10n6 & 5 & 6 & 1 & 4 & 0 & $\pi/10$ & 1.60849 & 1.93522 \\
10n7 & 3 & 7 & 3 & 4 & 0 & $\pi/6$ & 1.04300 & 0.52778 \\
10n8 & 3 & 5 & 4 & 13 & 0 & $\pi/6$ & 0.75398 & 2.81486 \\
10n9 & 3 & 8 & 2 & 7 & 0 & $\pi/6$ & 0.47752 & 1.45769 \\
10n10 & 3 & 8 & 2 & 5 & 0 & $\pi/6$ & 1.99805 & 2.06088 \\
10n11 & 4 & 7 & 3 & 8 & 0 & $\pi/8$ & 1.26920 & 0.79168 \\
10n12 & 3 & 8 & 2 & 7 & 0 & $\pi/6$ & 1.39486 & 2.85256 \\
10n13 & 3 & 5 & 2 & 11 & 0 & $\pi/6$ & 1.53309 & 1.28176 \\
10n14 & 3 & 5 & 2 & 9 & 0 & $\pi/6$ & 0.33929 & 1.58336 \\
10n15 & 3 & 8 & 2 & 3 & 0 & $\pi/6$ & 0.52778 & 1.58336 \\
10n16 & 3 & 8 & 6 & 7 & 0 & $\pi/6$ & 0.76654 & 0.52778 \\
10n17 & 3 & 10 & 1 & 4 & 0 & $\pi/6$ & 0.47752 & 1.70902 \\
10n18 & 3 & 4 & 1 & 10 & 0 & $\pi/6$ & 2.02318 & 1.91008 \\
10n19 & 3 & 5 & 2 & 7 & 0 & $\pi/6$ & 1.25663 & 2.04831 \\
10n20 & 3 & 4 & 5 & 10 & 0 & $\pi/6$ & 0.62831 & 1.93522 \\
10n21 & 3 & 5 & 2 & 3 & 0 & $\pi/6$ & 0.31415 & 1.58336 \\
10n22 & 3 & 10 & 4 & 13 & 0 & $\pi/6$ & 0.02513 & 0.02513 \\
10n23 & 3 & 7 & 8 & 11 & 0 & $\pi/6$ & 0.05026 & 1.28176 \\
10n24 & 3 & 4 & 7 & 14 & 0 & $\pi/6$ & 1.74672 & 1.48283 \\
10n25 & 3 & 7 & 1 & 12 & 0 & $\pi/6$ & 2.09858 & 1.58336 \\
10n26 & 3 & 7 & 5 & 6 & 0 & $\pi/6$ & 0.69115 & 1.58336 \\
10n27 & 4 & 5 & 1 & 5 & 0 & $\pi/8$ & 0.03769 & 4.80035 \\
10n28 & 3 & 4 & 7 & 11 & 0 & $\pi/6$ & 0.05026 & 2.23681 \\
10n29 & 3 & 5 & 2 & 11 & 0 & $\pi/6$ & 1.92265 & 2.09858 \\
10n30 & 3 & 7 & 2 & 4 & 0 & $\pi/6$ & 0.38955 & 1.52053 \\
10n31 & 3 & 7 & 4 & 15 & 0 & $\pi/6$ & 0.52778 & 0.75398 \\
10n32 & 3 & 5 & 7 & 10 & 0 & $\pi/6$ & 0.41469 & 2.31221 \\
10n33 & 3 & 5 & 2 & 9 & 0 & $\pi/6$ & 1.45769 & 1.58336 \\
10n34 & 3 & 8 & 1 & 12 & 0 & $\pi/6$ & 0.66601 & 1.58336 \\
10n35 & 3 & 4 & 10 & 13 & 0 & $\pi/6$ & 0.07539 & 3.00336 \\
10n36 & 4 & 5 & 1 & 4 & 0 & $\pi/8$ & 0.18849 & 1.18123 \\
10n37 & 3 & 8 & 1 & 4 & 0 & $\pi/6$ & 2.22424 & 0.22619 \\
10n38 & 3 & 4 & 5 & 10 & 0 & $\pi/6$ & 2.55097 & 1.49539 \\
10n39 & 3 & 5 & 3 & 10 & 0 & $\pi/6$ & 1.58336 & 0.52778 \\
10n40 & 3 & 5 & 7 & 14 & 0 & $\pi/6$ & 1.00530 & 2.55097 \\
10n41 & 3 & 10 & 2 & 7 & 0 & $\pi/6$ & 2.57610 & 2.78973 \\
10n42 & 3 & 10 & 1 & 8 & 0 & $\pi/6$ & 1.43256 & 0.79168 \\\hline
\end{tabular}
\end{center}
\caption{Fourier-$(1,1,2)$ descriptions of nonalternating non 2-bridge knots with 10 crossings. All amplitudes are 1.  Knot names are as in {\sl Knotscape}.}
\label{non 2-bridge 3 to 10 as fourier(1,1,2), part 4}
\end{table}%

\begin{table}[htdp]
\begin{center}
\begin{tabular}{|l|l|cccccccc|}
\hline
torus knot& knot&$n_x$ & $n_y$ & $n_{z,1}$&$n_{z,2}$ & $\phi_x$ & $\phi_y$ & $\phi_{z,1}$&$\phi_{z,2}$\\
\hline
$T_{2,3}$&3a1 & 2 & 3 & 1 & 2 & 0 & $\pi/4$ & $\pi/2$ & $\pi/4$ \\
$T_{2,5}$&5a2  & 2 & 5 & 2 & 3 & 0 & $\pi/4$ &$\pi/2$ & $\pi/4$ \\
$T_{2,7}$&7a7  & 2 & 7 & 2 & 5 & 0 & $\pi/4$ &$\pi/2$ & $\pi/4$ \\
$T_{2,9}$&9a41  & 2 & 9 & 2 & 7 & 0 & $\pi/4$ &$\pi/2$ & $\pi/4$ \\
$T_{2,11}$&11a367&2& 11& 2& 9& 0& $\pi/4$& $\pi/2$ & $\pi/4$ \\
$T_{2,13}$&13a4878&2& 13& 2& 11& 0& $\pi/4$& $\pi/2$ & $\pi/4$ \\
$T_{2,15}$&15a85263&2& 15& 2& 13& 0& $\pi/4$& $\pi/2$ & $\pi/4$ \\
$T_{3,4}$&8n3 & 3 & 4 & 1 & 3 & 0 & $\pi/6$ & 0.26389 & 1.58336 \\
$T_{3,5}$&10n21 & 3 & 5 & 2 & 3 & 0 & $\pi/6$ & 0.31415 & 1.58336 \\
$T_{3,7}$&14n21881 & 3 & 7 & 3 & 4 & 0 & $\pi/6$ & 1.57079 & 0.37699 \\
$T_{4,5}$&15n41185 & 4 & 5 & 1 & 4 & 0 & $\pi/8$ & 0.40212 & 1.58336 \\
$T_{3,8}$&16n783154 & 3 & 8 & 3 & 5 & 0 & $\pi/6$ & 1.57079 & 0.40212 \\

\hline
\end{tabular}
\end{center}
\caption{Fourier-$(1,1,2)$ descriptions of all torus knots up to 16 crossings. All amplitudes are 1.  Knot names are as in {\sl Knotscape}.}
\label{torus knots}
\end{table}%


\begin{thebibliography}{10}

\bibitem{BHJ1994}
M.~G.~V. Bogle, J.~E. Hearst, V.~F.~R. Jones, and L.~Stoilov.
\newblock Lissajous knots.
\newblock {\em J. Knot Theory Ramifications}, 3(2):121--140, 1994.

\bibitem{BZ:2003}
Gerhard Burde and Heiner Zieschang.
\newblock {\em Knots}, volume~5 of {\em de Gruyter Studies in Mathematics}.
\newblock Walter de Gruyter \& Co., Berlin, 2003.

\bibitem{C:2004}
Peter Cromwell.
\newblock {\em Knots and Links}.
\newblock Cambridge University Press, 2004.

\bibitem{HK:1979}
Richard Hartley and Akio Kawauchi.
\newblock Polynomials of amphicheiral knots.
\newblock {\em Math. Ann.}, 243(1):63--70, 1979.

\bibitem{HT:1998}
Jim Hoste and Morwen Thistlethwaite.
\newblock Knotscape.
\newblock {\em \tt http://www.math.utk.edu/$\sim$morwen}, 1998.

\bibitem{HTW:1998}
Jim Hoste, Morwen Thistlethwaite, and Jeff Weeks.
\newblock The first 1,701,936 knots.
\newblock {\em Math. Intelligencer}, 20(4):33--48, 1998.

\bibitem{HZ2006}
Jim Hoste and Laura Zirbel.
\newblock Lissajous knots with lissajous projections.
\newblock {\em arXiv: math.GT/0605632}, 2006.

\bibitem{JP1998}
Vaughan F.~R. Jones and J\'{o}zef~H. Przytycki.
\newblock Lissajous knots and and billiard knots.
\newblock {\em Banach Center Publications}, (42):145--163, 1998.

\bibitem{K}
Louis H.  Kauffman.
\newblock{Fourier Knots.}
\newblock{\em arXiv: q-alg/9711013}

\bibitem{Lamm}
Christoph Lamm.
\newblock Fourier knots.
\newblock {\em Preprint}.

\bibitem{Lamm1996}
Christoph Lamm.
\newblock There are infinitely many lissajous knots.
\newblock {\em Manuscripta Math.}, 93:29--37, 1996.

\bibitem{M:1971}
Kunio Murasugi.
\newblock On periodic knots.
\newblock {\em Comment. Math. Helv.}, 46:162--174, 1971.

\bibitem{P}
Jozef Przytycki.
\newblock Symmetric knots and billiard knots.
\newblock{\em arXiv: math.GT/0405151}

\bibitem{Rolfsen}
Dale Rolfsen.
\newblock {\em Knots and links}, volume~7 of {\em Mathematics Lecture Series}.
\newblock Publish or Perish Inc., Houston, TX, 1990.
\newblock Corrected reprint of the 1976 original.

\end{thebibliography}
\end{document}